\pdfoutput=1
\documentclass[12pt]{article}

\usepackage{mathptmx}       % selects Times Roman as basic font
\usepackage{helvet}         % selects Helvetica as sans-serif font
\usepackage{courier}        % selects Courier as typewriter font
\usepackage{type1cm}        % activate if the above 3 fonts are
\usepackage{makeidx}         % allows index generation
\usepackage{graphicx}        % standard LaTeX graphics tool
\usepackage{placeins}
\usepackage{multicol}        % used for the two-column index
\usepackage[bottom]{footmisc}% places footnotes at page bottom

\usepackage[UKenglish]{babel}

\usepackage{caption}
\usepackage{graphicx}
\usepackage{epstopdf}
\usepackage{todonotes}
\usepackage{listings}

\usepackage{amsmath}				% Basic AMS package
\usepackage{amssymb}				% Lots of new smbols
\usepackage{amsfonts}				% For blackboard bold, etc
\usepackage{amsthm}					% Better theorem environments
\usepackage{amscd}					% AMS Diagrams
\usepackage{mathtools}				% Extensions to amsmath package
\mathtoolsset{showonlyrefs=true}	% Only number equations referenced

\renewcommand{\vec}[1]{\mathbf{#1}} % Vectors as bold letters 

\newtheorem{definition}{Definition}
\newtheorem{proposition}{Proposition}

\newtheorem{theorem}{Theorem}

\usepackage{csquotes}
\usepackage[
    backend=biber,
    natbib=true,
    style=numeric-comp,
    sorting=none,
    maxnames=2,
    firstinits=true,
    url=true,
    arxiv=pdf,
    doi=true,
    eprint=true,
    isbn=false,
    date=terse,
	block=ragged,
	maxbibnames=99
]{biblatex}
\usepackage{tikz}
\usetikzlibrary{fit, shapes, backgrounds, calc, positioning}

\makeatletter
\tikzset{circle split part fill/.style args={#1,#2}{
 alias=tmp@name,
  postaction={
    insert path={
     \pgfextra{
     \pgfpointdiff{\pgfpointanchor{\pgf@node@name}{center}}
                  {\pgfpointanchor{\pgf@node@name}{east}}           
     \pgfmathsetmacro\insiderad{\pgf@x}
      \fill[#1] (\pgf@node@name.base) ([xshift=-\pgflinewidth]\pgf@node@name.east) arc
                          (0:180:\insiderad-\pgflinewidth)--cycle;
      \fill[#2] (\pgf@node@name.base) ([xshift=\pgflinewidth]\pgf@node@name.west)  arc
                           (180:360:\insiderad-\pgflinewidth)--cycle;
         }}}}}  
\makeatother

\tikzset{nix/.style={
	rectangle,
	inner sep=2pt,
    draw=none,
	fill=none
}}

\tikzset{Basic/.style={
	inner sep=0pt,
	minimum size=0.35cm,
	line width=0.3mm
}}

\tikzset{T/.style={
	shape=circle,
    Basic,
	draw=black,
	fill=black
}}

\tikzset{GT/.style={
	shape=circle,
    Basic,
	draw=gray,
	fill=gray
}}

\tikzset{DT/.style={
    shape=circle split,
	Basic,
	draw=black,
	rotate=-45
}}

\tikzset{GDT/.style={
    shape=circle split,
	Basic,
	draw=gray,
	rotate=45
}}

\tikzset{LT/.style={ T, circle split part fill={black,white}, rotate=45 } }
\tikzset{LLT/.style={ T, circle split part fill={black,white}, rotate=135 }}
\tikzset{RT/.style={ T, circle split part fill={black,white}, rotate=-45 } }
\tikzset{RRT/.style={ T, circle split part fill={black,white}, rotate=-135 }}

\tikzset{GLT/.style={ GT, circle split part fill={gray,white}, rotate=45 } }
\tikzset{GRT/.style={ GT, circle split part fill={gray,white}, rotate=-45 } }

\begin{document}

\title{A Randomized Tensor Train Singular Value Decomposition}
%\titlerunning{A Randomized Hierarchical SVD}
\author{Benjamin Huber, Reinhold Schneider and Sebastian Wolf}
%\authorrunning{B. Huber, R. Schneider and S. Wolf}

\maketitle

\abstract{The hierarchical SVD provides a quasi-best low rank approximation of high dimensional data in the hierarchical Tucker framework. Similar to the SVD for matrices, it provides a fundamental but expensive tool for tensor computations. In the present work we examine generalizations of randomized matrix decomposition methods to higher order tensors in the framework of the hierarchical tensors representation. In particular we present and analyze a randomized algorithm for the calculation of the hierarchical SVD (HSVD) for the tensor train (TT) format.
}

\section{Introduction}
Low rank matrix decompositions, such as the singular value decomposition (SVD) and the QR decomposition are principal tools in data analysis and scientific computing. For matrices with small rank both decompositions offer a tremendous reduction in computational complexity and can expose the underlying problem structure. In recent years generalizations of these low rank decompositions to higher order tensors have proven to be very useful and efficient techniques as well. In particular the hierarchical Tucker \cite{hackbusch2009new} and the tensor train \cite{oseledets2011tensor} format made quite an impact, as both formats allow to circumvent the notorious \emph{curse of dimensionality}, i.e. the exponential scaling of the ambient spaces with respect to the order of the tensors. Applications of these formats are as various as high-dimensional PDE's like the Fokker Planck equations and the many particle Schr\"odinger equations, applications in neuroscience, graph analysis, signal processing, computer vision and computational finance, see the extensive survey of \citet{grasedyck2013literature}. Also in a recent paper in machine learning \citet{cohen2015expressive} showed a connection between these tensor formats and deep neural networks and used this to explain the much higher power of expressiveness of deep neural networks over shallow ones.\\

One of the main challenges when working with these formats is the calculation of low rank decompositions of implicitly or explicitly given tensors, i.e.\ the high dimension analog of the classical SVD calculation. For matrices there exists a wide range of methods, which allow these calculations with high efficiently and precision. One particular branch are randomized methods which appear often in the literature, mostly as efficient heuristics to calculate approximate decompositions. It was only recently that thanks to new results form \emph{random matrix theory}, a rigorous analysis of these procedures became possible, see \cite{halko2011finding}. In this work we aim to extend some of these results for randomized matrix decompositions to the high dimensional tensor case. To this end we present an algorithm which allows the efficient calculation of the tensor train SVD (TT-SVD) for general higher order tensors. Especially for sparse tensors this algorithms exhibits a superior complexity, scaling only linear in the order, compared to the exponential scaling of the naive approach. Extending the results of \cite{halko2011finding}, we show that stochastic error bounds can also be obtained for these higher order methods.\\

This work focuses on the theoretical and algorithmic aspects of this randomized (TT-)SVD. However a particular application on our mind is the work in \cite{rauhut2015tensor, rauhut2016low}, where we treat the tensor completion problem. That is, in analogy to matrix completion, see e.g.\ \cite{candes2009exact, recht2011simpler, cai2010singular}, we want to reconstruct a tensor from $N$ measurements using a low rank assumption. We use an iterative (hard) thresholding procedure, which requires the (approximate) calculation of a low rank decomposition in each iteration of the algorithm. As the deterministic TT-SVD is already a fundamental tool, there are of course many further possible application for our randomized variant, see for example \cite{bachmayr2015adaptive, bachmayr2016tensor, bachmayr2016iterative}.\\

We start with a brief recap of tensor product spaces and introduce the notation used in the remainder of this work. In section \ref{sec:TensorDecomp} we give an overview of different tensor decompositions, generalizing the singular value decomposition from matrices to higher order tensors. In the second part a detailed introduction of the tensor train format is provided. Section \ref{sec:randMatrixDecomp} summarizes results for randomized matrix decompositions which are important for this work. In section \ref{sec:randTensorDecomp} we introduce our randomized TT-SVD scheme and prove stochastic error bounds for this procedure. An interesting relation between the proposed algorithm and the popular alternating least squares (ALS) algorithm is examined in section \ref{sec:als}. Section \ref{sec:numerics} collects several numerical experiments showing the performance of the proposed algorithms. Section \ref{sec:conclusions} closes with some concluding remarks.

\subsection{Tensor product spaces}
Let us begin with some preliminaries on tensors and tensor spaces. For an exhaustive introduction we refer to the monograph of \citet{hackbusch2012tensorBuch}.\\

Given Hilbert spaces $V_1, \ldots, V_d$ the tensor product space of order $d$
$$ \mathcal{V} = \bigotimes_{i=1}^d V_i \ , $$
is defined as the closure of the span of all elementary tensor products of vectors from $V_i$, i.e.\ 
$$
\mathcal{V} := \overline{ \mbox{span }  \left\{ \vec v_1 \otimes \vec v_2 \otimes \ldots \otimes \vec v_d \mid \vec v_i \in V_i \right\} } \ .
$$
The elements $\vec x \in \mathcal{V}$ are called tensors of order $d$. If each space $V_i$ is supplied with an orthonormal basis $\{ \varphi^i_{\mu_i} : \mu_i \in \mathbb{N} \}$, then any $\mathbf{x} \in \mathcal{V}$ can be represented as
$$
\vec x = \sum_{\mu_1=1 }^{\infty } \ldots \sum_{\mu_d = 1}^{\infty } \vec x [\mu_1 , \ldots , \mu_d] \  \varphi^1_{\mu_1}  \otimes \cdots \otimes  \varphi^d_{\mu_d} \ .
$$
Using this basis, with a slight abuse of notation, we can identify $\vec x \in \mathcal{V}$ with its representation by a $d$-variate function, often called hyper matrix, 
$$ 
\boldsymbol{\mu} = (\mu_1, \ldots , \mu_d ) \mapsto \vec x [\mu_1, \ldots , \mu_d] \in \mathbb{K} \ ,
$$
depending on discrete variables, usually called indices $\mu_i \in \mathbb{N}$. Analogous to vectors and matrices we use square brackets $\vec x[\mu_1, \ldots , \mu_d]$ to index the entries of this hypermatrix. Of course, the actual representation of $\vec x \in \mathcal{V}$ depends on the chosen bases $\varphi^i$ of $V_i$. The index $\mu_i$ is said to correspond to the $\mu$-th mode or equivalently the $\mu$-th dimension of the tensor. \\

In the remainder of this article, we confine to finite dimensional real linear spaces $V_i : = \mathbb{R}^{n_i}$, however most parts are easy to extend to the complex case as well. For these, the tensor product space
$$
\mathcal{V} = \bigotimes_{i=1}^d \mathbb{R}^{n_i} = \mathbb{R}^{n_1 \times n_2 \times \ldots \times n_d} :=  \mbox{span }  \left\{ \vec v_1 \otimes \vec v_2 \otimes \ldots \otimes \vec v_d \mid \vec v_i \in \mathbb{R}^{n_i} \right\}
$$
is easily defined. If it is not stated explicitly, the $V_i = \mathbb{R}^{n_i}$ are supplied with the canonical basis $ \{ \mathbf{e}^i_1 , \ldots, \mathbf{e}^i_{n_i} \} $ of the vector spaces $\mathbb{R}^{n_i}$. Then every $\vec x \in \mathcal{V}$ can be represented as
\begin{equation}
\vec x = \sum_{\mu_1=1 }^{n_1} \ldots \sum_{\mu_d = 1}^{n_d} \vec x [ \mu_1 , \ldots , \mu_d ] \;
 \mathbf{e}^1_{\mu_1} \otimes \cdots \otimes \mathbf{e}^d_{\mu_d} \ . \label{eq:elem}
\end{equation}
We equip the finite dimensional linear space $\mathcal{V}$ with the inner product 
$$
\langle \vec x , \vec y \rangle := \sum_{\mu_1= 1}^{n_1}  \cdots \sum_{\mu_d =1}^{n_d} \vec x [\mu_1, \ldots , \mu_d ] \ \vec y [\mu_1, \ldots , \mu_d ] \ . 
$$
and the corresponding $l_2$-norm $\| \vec x  \| = \sqrt{\langle \vec x, \vec x \rangle}$. We use the fact that for a Hilbert space $V_i$ the dual space $V^*$ is isomorphic to $V_i$ and use the identification $V \simeq V^*$. For the treatment of reflexive Banach spaces we refer to \cite{falco2012minimal, falco2015geometric}.\\

The number of possibly non-zero entries in the representation of $\vec x$ is $n_1 \cdots n_d = \Pi_{i=1}^d n_i$, and with $n= \max \{n_i : i=1, \ldots, d \}$, the dimension of the space $\mathcal{V}$ scales exponentially in $d$, i.e. $ \mathcal{O} (n^d) $. This is often referred to as the {\em curse of dimensions} and presents the main challenge when working with higher order tensors.

\subsection{Tensor Contractions and Diagrammatic Notation}
Important concepts for the definitions of tensor decompositions are so called \emph{matricizations} and \emph{contractions} introduced in this section. \\

The matricization or flattening of a tensor is the reinterpretation of the given tensor as a matrix, by combining a subset of modes to a single mode and combining the remaining modes to a second one. 
\begin{definition}[Matricization or Flattening]
	Let $[n]=\{1,2,\ldots,n\}$ and $\vec x \in \mathbb{R}^{n_1 \times \ldots \times n_d}$ be a tensor of order $d$. Furthermore let $\alpha \subseteq [d]$ be a subset of the modes of $\vec x$, and let $\beta = [d] \backslash \alpha$ be its complement. Given two bijective functions $\mu_\alpha : [n_{\alpha_1}]\times [n_{\alpha_2}]\times\ldots\rightarrow [n_{\alpha_1} \cdot n_{\alpha_2}\cdots]$ and $\mu_\beta$ respectively.
    
    The $\alpha$-matricization or $\alpha$-flattening 
	\begin{align}
		\hat M_\alpha :  \mathbb{R}^{n_1 \times  \ldots \times n_d} \rightarrow  \mathbb{R}^{m_\alpha \times m_\beta}\\
		\vec x \mapsto \hat M_\alpha(\vec x)
	\end{align}
	of $\vec x$ is defined entry-wise as
	\begin{align}
		\vec x[i_1, \ldots, i_d] &=: \hat M_\alpha(\vec x)\left[ \mu_\alpha(i_{\alpha_1}, i_{\alpha_2} \ldots), \mu_\beta(i_{\beta_1},i_{\beta_2},\ldots) \right] .
	\end{align}
    A common choice for $\mu_\alpha$ and $\mu_\beta$ is $\mu(i_1, i_2, \ldots) = \sum_k i_k \prod_{j>k} n_j$. The actual choice is of no significance though, as long as it stays consistent.
    The matrix dimensions are given as $m_\alpha = \prod_{j \in \alpha} n_j$ and $m_\beta = \prod_{j \in \beta} n_j$.
\end{definition}

The inverse operation is the de-matricization or unflattening $\hat M^{-1}$. In principle it is possible to define de-matricization for any kind of matrix, typically called tensorization. However this requires to give the dimensions of the resulting tensor and the details of the mapping alongside with the operator. Instead, in this work the de-matricization is only applied to matrices where at least one mode of the matrix encodes a tensor structure through a former matricization, in which the details of the mapping are clear from the context. For all other modes the de-matricization is simply defined to be the identity.\\

The second important tool are tensor contractions, which are generalizations of the matrix-vector and matrix-matrix multiplications to higher order tensors.
\begin{definition}[Tensor Contraction]
Let $\vec x \in \mathbb{R}^{n_1 \times \ldots\times n_d}$ and $\vec y \in \mathbb{R}^{m_1 \times \ldots \times m_e}$ be two tensors of order $d$ and $e$ respectively, with $n_k = m_l$. The contraction of the $k$-th mode of $\vec u$ with the $l$-th mode of $\vec v$
\begin{align}
\vec z := \vec x \circ_{k,l} \vec y
\end{align}
is defined entry-wise as
\begin{align*}
& \ \ \ \ \ \ \vec z[i_1, \ldots, i_{k-1}, i_{k+1}, \ldots, i_d, j_1, \ldots, j_{l-1},j_{l+1}, \ldots, j_e]\\
&= \sum_{p=0}^{n_k} \vec x[i_1, \ldots, i_{k-1}, p, i_{k+1}, \ldots, i_d]  \ \vec y[ j_1, \ldots, j_{l-1},p , j_{l+1}, \ldots, j_e]
\end{align*}
or via the matricizations
$$ \vec z = \hat M^{-1} \left( \hat M_{\{k\}} (\vec x)^T \hat M_{\{l\}}(\vec y) \right) \ .$$
The resulting tensor $\vec z \in \mathbb{R}^{{n_1}  \times \ldots \times n_{k-1}  \times n_{k+1} \times \ldots  \times n_d \times m_1 \times \ldots m_{k-1} \times m_{k+1} \times \ldots \times m_{e}}$ is of order $d+e-2$. Note that in order for this operation to be well-defined, $n_k = m_l$ must hold.
\end{definition}

If no indices are specified, i.e.\ only $\circ$, a contraction of the last mode of the left operand and the first mode of the right operand is assumed. If tuples of indices are given, e.g.\ $\circ_{(i,j,k), (l,p,q)}$, a contraction of the respective mode pairs ($i/l$, $j/p$ $k/q$) is assumed.\footnote{As one can easily show the order of the contractions does not matter.} As writing this for larger tensor expressions quickly becomes cumbersome we also use a diagrammatic notation to visualize the contraction. In this notation a tensor is depicted as a dot or box with edges corresponding to each of its modes. If appropriate the cardinality of the corresponding index set is given as well. From left to right the following shows this for an order one tensor (vector) $\vec v \in \mathbb{R}^n$, an order two tensor (matrix) $\vec A \in \mathbb{R}^{m \times n}$ and an order 4 tensor $\vec x \in \mathbb{R}^{n_1 \times n_2 \times n_3  \times n_4 }$.
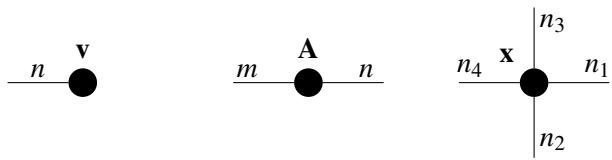
\begin{center}
	\begin{tikzpicture}
		\node[T, label=above:{$\vec v$}] (p1) at (0,0){};
		\draw (p1) -- node[nix, above]{$n$} ++(-1,0);

		\node[T, label=above:{$\vec A$}] (p2) at (3,0){};
		\draw (p2) -- node[nix, above left]{$m$} ++(-1,0);
		\draw (p2) -- node[nix, above right]{$n$} ++(1,0);

		\node[T, label=above left:{$\vec x$}] (p3) at (6,0){};
		\draw (p3) -- node[nix, above left]{$n_4$} ++(-1,0);
		\draw (p3) -- node[nix, above right]{$n_1$} ++(1,0);
		\draw (p3) -- node[nix, below right]{$n_2$} ++(0,-1);
		\draw (p3) -- node[nix, above right]{$n_3$} ++(0,1);
	\end{tikzpicture}
\end{center}
If a contraction is performed between the modes of two tensors the corresponding edges are joined. The following shows this exemplary for the inner product of two vectors $\vec u, \vec v \in \mathbb{R}^n$ and a matrix-vector product with $\vec A \in \mathbb{R}^{m \times n}$ and $\vec v \in  \mathbb{R}^n$.
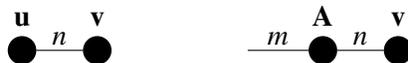
\begin{center}
	\begin{tikzpicture}
		\node[T, label=above:{$\vec u$}] (p1) at (0,0){};
		\node[T, label=above:{$\vec v$}] (p2) at (1,0){};

		\node[T, label=above:{$\vec A$}] (p3) at (4,0){};
		\draw (p3) -- node[nix, above]{$m$} ++(-1,0);
		\node[T, label=above:{$\vec v$}] (p4) at (5,0){};

		\path 
		(p1) edge node[nix, above]{$n$}(p2)
		(p3) edge node[nix, above]{$n$}(p4);
	\end{tikzpicture}
\end{center}

There are two special cases concerning orthogonal and diagonal matrices. If a specific matricization of a tensor yields an orthogonal or diagonal matrix, the tensor is depicted by a half filled circle (orthogonal) or a circle with a diagonal bar (diagonal) respectively. The half filling and the diagonal bar both divide the circle in two halves. The edges joined to either half, correspond to the mode sets of the matricization, which yields the orthogonal or diagonal matrix. As an example the diagrammatic notation can be used to depict the singular value decomposition $\vec A = \vec U \vec S \vec V^T$ of a matrix $\vec A \in \mathbb{R}^{m \times n}$ with rank $r$, as shown in the following.
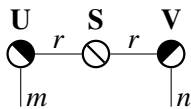
\begin{center}
 \begin{tikzpicture}
	\node[RT, label=above left:{$\vec U$}] (u) at (1,0){};
	\draw (u) -- node[nix, below right]{$m$} ++(0,-0.7);
	\node[DT, label=above left:{$\vec S $}] (s) at (2,0){};
	\node[LT, label=above right:{$\vec V$}] (v) at (3,0){};
	\draw (v) -- node[nix, below right]{$n$} ++(0,-0.7);

	\path 
	(u) edge node[nix, above]{$r$}(s)
	(s) edge node[nix, above]{$r$}(v);
\end{tikzpicture}
\end{center}

\section{Low Rank Tensor Decompositions}
\label{sec:TensorDecomp}

In this section we give an introduction to the low rank tensor decomposition techniques used in the remainder of this work. As there are in fact quite different approaches to generalize the singular value decomposition, and thereby also the definition of the rank, to higher order tensors, we start with an overview of the most popular formats. For an in depth overview including application we refer to the survey of \citet{grasedyck2013literature}. In the second part we provide a detailed introduction of the tensor train format, which is used in the remainder of this work.\\

The probably best known and classical tensor decomposition is the representation by a sum of elementary tensor products, i.e.\ 
\begin{equation}
\vec x = \sum_{i=1}^r \vec u_{1, i} \otimes \vec u_{2, i} \otimes \ldots \otimes \vec u_{d, i} \label{eq:canonicalFormat}
\end{equation}
where $\vec x \in \mathbb{R}^{n_1 \times \ldots \times n_d}$ and $\vec u_{k, i} \in \mathbb{R}^{n_k}$ are vectors from the respective vector spaces. This format is mainly know as the \emph{canonical format} but also appears in the literature under the names canonical polyadic (CP) format, CANDECOMP and PARAFAC. The canonical or CP rank is defined as the minimal $r$ such that a decomposition as in \eqref{eq:canonicalFormat} exists. Note that in general there is no unique CP representation with minimal rank. This is somewhat expected, since even for matrices the SVD is not unique if two or more singular values coincide. Some discussion on the uniqueness can be found in the paper of \citet{kolda2009tensor}. For tensors with small canonical rank \eqref{eq:canonicalFormat} offers a very efficient representation, requiring only $\mathcal{O}(rdn)$ storage instead of $\mathcal{O}(n^d)$ for the direct representation. Unfortunately the canonical format suffers from several difficulties and instabilities. First of all the task of determining the canonical rank of a tensor with order $d>2$ is, in contrast to matrices, highly non trivial. In fact it was shown by \cite{haastad1990tensor} that even for order $d=3$, the problem of deciding whether a rational tensor has CP-rank $r$ is NP-hard (and NP-complete for finite fields). Consequently also the problem of calculating low rank approximations proves to be challenging. That is, given a tensor $\vec x \in \mathbb{R}^{n_1 \times \ldots \times n_d}$ and a CP-rank $r$, finding the best CP-rank $r$ approximation
\begin{align}
	\vec x^* = \underset{\vec y \in \mathbb{R}^{n_1 \times \ldots \times n_d}, \  \text{CP-rank}(\vec y) \leq r}{\text{argmin}} \left( \lVert \vec x - \vec y \rVert \right) \ . \label{eq:cpBestApp} 
\end{align}
The norm $\lVert \cdot \rVert $ used, may differ depending on the application. In the matrix case the Eckart-Young theorem provides that for the Frobenius and spectral norm this best approximation can be straightforwardly calculated by a truncated SVD. In contrast, \citet{keinCP} proved that the problem of the best CP-rank $r$ approximation, as formulated in \eqref{eq:cpBestApp}, is ill-posed for many ranks $r \geq 2$ and all orders $d > 2$ regardless of the choice of the norm $\lVert \cdot \rVert$. Furthermore they showed that the set of tensors that do not have a best CP-rank $r$ approximation is a non-null set, i.e.\ there is a strictly positive probability that a randomly chosen tensor does not admit a best CP-rank $r$ approximation. Finally it was shown by \citet{keinCP} that neither the set $\{ \vec x \in  \mathbb{R}^{n_1 \times  \ldots \times n_d} \ \mid \ \text{CP-rank}(\vec x) = r \}$ of all tensors with CP-rank $r$, nor the set $\{ \vec x \in \mathbb{R}^{n_1 \times \ldots \times n_d} \ \mid \ \text{CP-rank}(\vec x) \leq  r \} $ of all tensors with CP-rank at most $r$ is closed for $d > 2$. These are some severe difficulties for both the theoretical and practical work with the canonical format.\\

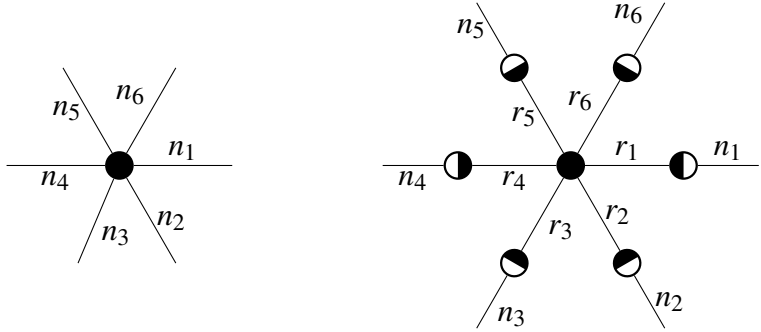
\begin{figure}
\centering	
\begin{tikzpicture}
		\node[T] (x1) at (-6,0){};
		\draw (x1) -- node[nix, above]{$n_1$} ++(1.5,0);
		\draw (x1) -- node[nix, right]{$n_2$} ++(0.75,-1.299);
		\draw (x1) -- node[nix, below right]{$n_3$} ++(-0.55,-1.299);
		\draw (x1) -- node[nix, below]{$n_4$} ++(-1.5,0);
		\draw (x1) -- node[nix, left]{$n_5$} ++(-0.75,1.299);
		\draw (x1) -- node[nix, above left]{$n_6$} ++(0.75,1.299);

		\node[T] (p1) at (0,0){};
		\node[T, circle split part fill={black,white}, rotate=90] (p2) at (1.5,0){};
		\draw (p2) -- node[nix, above]{$n_1$} ++(1,0);
		\node[T, circle split part fill={black,white}, rotate=30] (p3) at (0.75, -1.299){};
		\draw (p3) -- node[nix, right]{$n_2$} ++(0.5,-0.866);
		\node[T, circle split part fill={black,white}, rotate=330] (p4) at (-0.75, -1.299){};
		\draw (p4) -- node[nix, below right]{$n_3$} ++(-0.5,-0.866);
		\node[T, circle split part fill={black,white}, rotate=270] (p5) at (-1.5,0){};
		\draw (p5) -- node[nix, below]{$n_4$} ++(-1,0);
		\node[T, circle split part fill={black,white}, rotate=210] (p6) at (-0.75, 1.299){};
		\draw (p6) -- node[nix, left]{$n_5$} ++(-0.5,0.866);
		\node[T, circle split part fill={black,white}, rotate=150] (p7) at (0.75, 1.299){};
		\draw (p7) -- node[nix, above left]{$n_6$} ++(0.5,0.866);

		\path 
		(p1) edge node[nix, above]{$r_1$}(p2)
		(p1) edge node[nix, right]{$r_2$}(p3)
		(p1) edge node[nix, below right]{$r_3$}(p4)
		(p1) edge node[nix, below]{$r_4$}(p5)
		(p1) edge node[nix, left]{$r_5$}(p6)
		(p1) edge node[nix, above left]{$r_6$}(p7);

	 \end{tikzpicture}
\caption[Diagrammatic notation of the Tucker decomposition]{Left: A tensor $\vec x \in  \mathbb{R}^{n_1 \times  \ldots \times n_6}$ of order $6$. Right: Its Tucker decomposition.}
\label{fig:tucker}
\end{figure}

The second classical approach to generalize the SVD to higher order tensors is the subspace based Tucker decomposition. It was first introduced by \citet{tucker1966} in $1963$ and has been refined later on in many works, see e.g. \cite{kolda2009tensor, de2000multilinear, hackbusch2012tensorBuch}. Given a tensor $\vec x \in \mathbb{R}^{n_1 \times \ldots \times n_d}$, the main idea is to find minimal subspaces $U_i \subseteq \mathbb{R}^{n_i}$, such that $\vec x \in \mathbb{R}^{n_1 \times \ldots \times n_d}$ is an element of the induced tensor space
\begin{align}
	\mathcal{U} = \bigotimes _{i=1}^d U_i \subseteq \bigotimes_{i=1}^d \mathbb{R}^{n_i} = \mathbb{R}^{n_1 \times \ldots \times n_d} \ .
\end{align}
Let $r_i = \text{dim}(U_i)$ denote the dimension of the $i$-th subspace and let $\left\{\vec u_{i,j}, j=0, \ldots, r_i \right\}$ be an orthonormal basis of $U_i$. If the subspaces are chosen such that $\vec x \in \mathcal{U}$, then \eqref{eq:elem} states that there is a $\vec c$ such that $\vec x$ can be expressed as
\begin{align}
	\vec x = \sum_{\nu_1=1}^{r_1}  \ldots \sum_{\nu_d}^{r_d} \vec c[\nu_1, \nu_2, \ldots, \nu_d] \cdot \vec u_{1, \nu_1} \otimes \ldots \otimes \vec u_{d, \nu_d} \ .
\end{align}
Usually the basis vectors are combined to orthogonal matrices $\vec U_i = \left( \vec u_{i,1}, \ldots, \vec u_{i,r_i}\right)$, called \emph{basis matrices}. This leads to the following common form of the Tucker format
\begin{align}
		\vec x[\mu_1, \ldots, \mu_d] &= \sum_{\nu_1=1}^{r_1}  \ldots \sum_{\nu_d}^{r_d} \vec c[\nu_1, \ldots, \nu_d] \vec U_1[\mu_1, \nu_1] \ldots \vec U_d[\mu_d, \nu_d] \ . \label{eq:tuckerDecomp} 
\end{align}
The order $d$ tensor $\vec c \in \mathbb{R}^{r_1 \times \ldots \times r_d}$ of the prefactors is usually called \emph{core tensor}. The $d$-tuple $\vec r = (r_1, r_2, \ldots, r_d)$ of the subspace dimensions is called the representation rank and is associated with the particular representation. The \emph{Tucker rank} (T-rank) of $\vec x$ is defined as the unique minimal $d$-tuple $\vec r^* = (r_1^*, \ldots, r_d^*)$, such that there exists a Tucker representation of $\vec x$ with rank $\vec r^*$. Equation \eqref{eq:tuckerDecomp} consists of $d$ tensor contractions, that can be visualized in the diagrammatic notation, which is exemplarily shown in Figure \ref{fig:tucker} for $d=6$. Note that even for the minimal T-rank, the Tucker decomposition is not unique, as for any orthogonal matrix $\mathbf{Q}_i \in \mathbb{R}^{r_i \times r_i}$, one can define a matrix $\tilde{\mathbf{U}}_i = \mathbf{U}_i \mathbf{Q}_i$ and the tensor 
$$
\tilde{\mathbf{c}}\left[ \mu_1, \ldots, \mu_d \right] =\sum_{\nu=1}^{r_i  } \mathbf{c} \left[ \mu_1,\ldots,\nu,\ldots \mu_d \right] \mathbf{Q}^{T} \left[ \nu, \mu_i \right]
$$
such that the tensor $\vec x$ can also be written as
$$
\mathbf{u} \left[ \mu_1,\ldots, \mu_d \right] = \sum_{\nu_1=1}^{r_1} \ldots \sum_{v _d=1}^{r_d} \tilde{\mathbf{c}} \left[ \nu_1,\ldots, \nu_d \right] \vec U_1[\mu_1, \nu_1] \ldots \tilde{\vec U}_i[\mu_i, \nu_i] \ldots \vec U_d[\mu_d, \nu_d],
$$
which is a valid Tucker decomposition with the same rank.\\

It is shown by \citet{de2000multilinear} that the Tucker rank as the minimal $d$-tuple is indeed well-defined and that the entries $r_i$ of the Tucker rank correspond to the rank of the $i$-th mode matricization of the tensor. That is
\begin{align}
	\text{T-rank} (\vec x) = \left( \text{rank}(\hat M_{\{1\}} (\vec x)), \ldots, \text{rank}(\hat M_{\{d\}} (\vec x))\right) \ . \label{eq:tuckerRank} 
\end{align}
The proof is tightly linked to the fact, that a Tucker representation of a tensor $\vec x \in \mathbb{R}^{n_1 \times  \ldots \times n_d}$ with minimal representation rank, can be obtained by successive singular value decompositions. This procedure is referred to as the higher order singular value decomposition (HOSVD), see \cite{de2000multilinear} for the details. Using truncated SVDs an approximation of $\vec x$ by a tensor $\vec x^*$ with lower T-rank $\vec r^* = (r_1^*,  \ldots r_d^*) \preceq (r_1, \ldots, r_d)$, can be obtained. Where the symbol $\preceq$ denotes an entry-wise $\leq$, i.e.\  $(r_1, \ldots r_d) \preceq (r_1^*, \ldots, r_d^*) \iff r_i \leq r_i^* \forall i$. In contrast to the Eckart-Young theorem for matrices the approximation $\vec x^*$ obtained in this way is \emph{not} the best T-rank $\vec r^*$ approximation of $\vec x$. However it is a quasi-best approximation by a factor $\sqrt{d}$, i.e.
\begin{align}
	\lVert \vec x-\vec x^* \rVert_F \leq \sqrt{d} \min_{\vec y\ : \ \operatorname{T-rank}(\vec y) \preceq \vec r^*} \left( \lVert \vec x - \vec y \rVert_F \right) \ .
\end{align}
For many applications this quasi-best approximation is sufficient. As for the canonical format, finding the true best approximation is at the very least NP-hard in general, as it is shown by \cite{hillar2013most} that finding the best rank $(1, \ldots, 1)$ approximation is already NP-hard. To store a tensor in the Tucker format only the core tensor and the basis matrices have to be stored. This amounts to a storage requirement of $\mathcal{O} (r^d + dnr)$, where $r := \max(r_1, \ldots, r_d)$ and $n := \max(n_1, \ldots, n_d)$. Compared to the $\mathcal{O}(n^d)$ this is a major reduction but does not break the curse of dimensionality as the exponential scaling in $d$ remains.\\

A more recent development is the \emph{hierarchical Tucker (HT)} format, introduced by \citet{hackbusch2009new}. It inherits most of the advantages of the Tucker format, in particular a generalized higher order SVD, see \cite{grasedyck2010hierarchical}. But in contrast to the Tucker format the HT format allows a linear scaling with respect to the order for the storage requirements and common operations for tensors of fixed rank. The main idea of the HT format is to extend the subspace approach of the Tucker format by a multi-layer hierarchy of subspaces. For an in-depth introduction of the hierarchical Tucker format we refer to the pertinent literature, e.g.\ \cite{hackbusch2009new, grasedyck2010hierarchical, hackbusch2012tensorBuch}. In this work we will instead focus on the \emph{tensor train (TT)} format, as introduced by \citet{oseledets2011tensor}. The TT format offers mostly the same advantages as the more general HT format, while maintaining a powerful simplicity. In fact it can to some extend be seen as a special case of the HT format, see \citet{grasedyck2011introduction} for details on the relation between the TT and HT format.

\subsection{Tensor Train Format}
\label{sec:ttFormat}

In this subsection we give a detailed introduction to the \emph{tensor train (TT)} format. In the formulation used in this work the TT format was introduced by \citet{oseledets2011tensor}, however an equivalent formulation was known in quantum physic for quite some time, see e.g.\ \cite{perez2006matrix} for an overview. The idea of the TT-format is to separate the modes of a tensor into $d$ order two and three tensors. This results in a tensor network that is exemplary shown for an order four tensor $\vec x = \vec W_1 \circ \vec W_2 \circ \vec W_3 \circ \vec W_4$ in the following diagram.
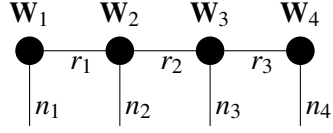
\begin{center}
	\begin{tikzpicture}
		\node[T, label=above:{$\vec W_1$}] (u1) at (0,0){};
		\node[T, label=above:{$\vec W_2$}] (u2) at (1.2,0){};
		\node[T, label=above:{$\vec W_3$}] (u3) at (2.4,0){};
		\node[T, label=above:{$\vec W_4$}] (u4) at (3.6,0){};
		\draw (u1) -- node[nix, below right]{$n_1$} ++(0,-1);
		\draw (u1) -- node[nix, below]{$r_1$} ++(1.2, 0);
		\draw (u2) -- node[nix, below right]{$n_2$} ++(0, -1);
		\draw (u2) -- node[nix, below]{$r_2$} ++(1.2, 0);
		\draw (u3) -- node[nix, below right]{$n_3$} ++(0,-1);
		\draw (u3) -- node[nix, below]{$r_3$} ++(1.2, 0);
		\draw (u4) -- node[nix, below right]{$n_4$} ++(0, -1);
	 \end{tikzpicture}
\end{center}
Formally it can be defined as follows.
\begin{definition}[Tensor Train Format] \label{def:TT} 
	Let $\vec x \in \mathbb{R}^{n_1 \times \ldots \times n_d}$ be a tensor of order $d$. A factorization
	\begin{align}
		\vec x = \vec W_1 \circ \vec W_2 \circ \ldots \circ \vec W_{d-1} \circ \vec W_d \ , \label{eq:TTRepresentation} 
	\end{align} 
	of $\vec x$, into component tensors $\vec W_1 \in \mathbb{R}^{n_1 \times r_1}$, $\vec W_i \in \mathbb{R}^{r_{i-1} \times n_i \times  r_i}$ ($i=2, \ldots, d-1$) and $\vec  W_d \in \mathbb{R}^{r_{d-1} \times n_d}$, is called a tensor train (TT) representation of $\vec x$. Equivalently \eqref{eq:TTRepresentation} can be given entry-wise as
	\begin{align*}
		&\vec x[i_1, \ldots, i_d] = \\
        &\sum_{j_1} \ldots \sum_{j_{d-1}} \vec W_1[i_1, j_1] \ \vec W_2[j_1,i_2,j_2] \ldots \vec W_{d-1}[j_{d-2}, i_{d-1}, j_{d-1}] \ \vec W_d[j_{d-1}, i_d] \ .
	\end{align*}
	The tuple of the dimensions $\vec r = (r_1, \ldots, r_{d-1})$ of the component tensors is called the representation rank and is associated with the specific representation. In contrast the tensor train rank (TT-rank) of $\vec x$ is defined as the minimal rank tuple $\vec r^* = (r_1^*, \ldots, r_{d-1}^*)$ such that there exists a TT representation of $\vec x$ with representation rank equal to $\vec r^*$.
\end{definition}

As for the Tucker format the TT-rank is well defined and linked to the rank of specific matricizations via 
$$
\text{TT-Rank}(\vec x) = \left( \text{rank}(\hat M_{\{1\}}(\vec x)), \text{rank}(\hat M_{\{1, 2\}}(\vec x)), \ldots, \text{rank}(\hat M_{\{1, 2, \ldots, d-1\}}(\vec x)) \right) \ .
$$
The proof is again closely linked to the fact that a tensor train decomposition of an arbitrary tensor can be calculated using successive singular value decompositions. This procedure is commonly referred to as the TT-SVD. For this work the TT-SVD is of particular importance as it constitutes the deterministic baseline for our randomized approach in section \ref{sec:randTensorDecomp}. In the following we therefore provide a complete step by step description of this procedure.\\

\begin{figure}
\centering
\begin{tikzpicture}
	\pgfmathsetmacro{\spacing}{3.1};
	\pgfmathsetmacro{\Nshift}{1.3};
	\pgfmathsetmacro{\shift}{1.15};
	\pgfmathsetmacro{\posA}{1};
	\pgfmathsetmacro{\posB}{6.2};
	\pgfmathsetmacro{\posC}{11.6-\Nshift};
	\pgfmathsetmacro{\posD}{1-0.6*\Nshift};
	\pgfmathsetmacro{\posE}{6.2-0.7*\Nshift};
	\pgfmathsetmacro{\posF}{11.6-\Nshift-0.5};
	\pgfmathsetmacro{\posG}{3-2*\Nshift};
	\pgfmathsetmacro{\posH}{10.6-1.5*\Nshift-0.4};
	\pgfmathsetmacro{\posI}{3-\Nshift-0.5};
	\pgfmathsetmacro{\posJ}{10.6-1.5*\Nshift-0.5};
	\pgfmathsetmacro{\posK}{6.5-2.5*\Nshift};
	\pgfmathsetmacro{\posL}{3-2*\Nshift};
	\pgfmathsetmacro{\posM}{11-1.5*\Nshift};

	\node[T, label=above right:{$\vec x_0$}] (t1) at (\posA,0){};
	\draw (t1) -- node[nix, below right]{$n_3$} ++(1,0);
	\draw (t1) -- node[nix, below right]{$n_4$} ++(0,-1);
	\draw (t1) -- node[nix, below left]{$n_1$} ++(-1,0);
	\draw (t1) -- node[nix, above left]{$n_2$} ++(0,1);
	\node at (\posA, -1.4) {(a)};

	\node[T, label=above:{$\hat M_{(1)}(\vec x_0)$}] (t) at (\posB,0){};
	\draw (t) -- node[nix, below left]{$n_1$} ++(-1,0);
	\draw (t) -- node[nix, below right]{$n_2 n_3 n_4$} ++(1, 0);
	\node at (\posB, -1.4) {(b)};

	\node[RT, label=above left:{$\vec U_1$}] (u1) at (\posC,0){};
	\node[DT, label=above left:{$\vec S _1$}] (s) at (\posC+\Nshift,0){};
	\node[LT, label=above right:{$\vec V_1$}] (v) at (\posC+2*\Nshift,0){};
	\draw (u1) -- node[nix, below right]{$n_1$} ++(0,-1);
	\draw (u1) -- node[nix, below]{$r_1$} ++(\shift, 0);
	\draw (s) -- node[nix, below]{$r_1$} ++(\shift, 0);
	\draw (v) -- node[nix, below right]{$n_2n_3n_4$} ++(0, -1);
	\node at (\posC+\Nshift, -1.4) {(c)};

	\node[RT, label=above left:{$\vec W_1$}] (u1) at (\posD,-\spacing){};
	\node[T, label=above:{$\hat M_{(1)}(\vec x_1)$}] (t) at (\posD+\Nshift,-\spacing){};
	\draw (u1) -- node[nix, below right]{$n_1$} ++(0,-1);
	\draw (u1) -- node[nix, below]{$r_1$} ++(\shift, 0);
	\draw (t) -- node[nix, below right]{$n_2n_3n_4$} ++(0, -1);
	\node at (\posD +0.6*\Nshift, -\spacing-1.4) {(d)};

	\node[RT, label=above left:{$\vec W_1$}] (u1) at (\posE,-\spacing){};
	\node[T, label=above right:{$\vec x_1$}] (t) at (\posE+\Nshift,-\spacing){};
	\draw (u1) -- node[nix, below right]{$n_1$} ++(0,-1);
	\draw (u1) -- node[nix, below]{$r_1$} ++(\shift, 0);
	\draw (t) -- node[nix, above left]{$n_2$} ++(0, 1);
	\draw (t) -- node[nix, below right]{$n_3$} ++(1, 0);
	\draw (t) -- node[nix, below right]{$n_4$} ++(0, -1);
	\node at (\posE +0.7*\Nshift, -\spacing-1.4) {(e)};

	\node[GRT, label=above left:{$\vec W_1$}] (u1) at (\posF,-\spacing){};
	\node[T, label=above:{$\hat M_{(1,2)}(\vec x_1)$}] (t) at (\posF+2*\Nshift,-\spacing){};
	\draw (u1) -- node[nix, below right]{$n_1$} ++(0,-1);
	\draw (u1) -- node[nix, below]{$r_1$} ++(0.8*\shift, 0);
	\draw (t) -- node[nix, below]{$r_1n_2$} ++(-\shift, 0);
	\draw (t) -- node[nix, below right]{$n_3n_4$} ++(1, 0);
	\node at (\posF +\Nshift+0.5, -\spacing-1.4) {(f)};

	\node[GRT, label=above left:{$\vec W_1$}] (u1) at (\posG,-2*\spacing){};
	\node[RT, label=above left:{$\vec U_2$}] (u2) at (\posG+2*\Nshift,-2*\spacing){};
	\node[DT, label=above left:{$\vec S_2$}] (s) at (\posG+3*\Nshift,-2*\spacing){};
	\node[LT, label=above right:{$\vec V_2$}] (v) at (\posG+4*\Nshift,-2*\spacing){};
	\draw (u1) -- node[nix, below right]{$n_1$} ++(0,-1);
	\draw (u1) -- node[nix, below]{$r_1$} ++(0.8*\shift, 0);
	\draw (u2) -- node[nix, below]{$r_1n_2$} ++(-1*\shift, 0);
	\draw (u2) -- node[nix, below]{$r_2$} ++(\shift, 0);
	\draw (s) -- node[nix, below]{$r_2$} ++(\shift, 0);
	\draw (v) -- node[nix, below right]{$n_3n_4$} ++(0, -1);
	\node at (\posG +2*\Nshift, -2*\spacing-1.4) {(g)};

	\node[GRT, label=above left:{$\vec W_1$}] (u1) at (\posH,-2*\spacing){};
	\node[RT, label=above left:{$\hat M_{(1,2)}(\vec W_2)$}] (u2) at (\posH+2*\Nshift,-2*\spacing){};
	\node[T, label=above:{$\hat M_{(1)}(\vec x_2)$}] (t) at (\posH+3*\Nshift+0.8,-2*\spacing){};
	\draw (u1) -- node[nix, below right]{$n_1$} ++(0,-1);
	\draw (u1) -- node[nix, below]{$r_1$} ++(0.8*\shift, 0);
	\draw (u2) -- node[nix, below]{$r_1n_2$} ++(-1*\shift, 0);
	\draw (u2) -- node[nix, below]{$r_2$} ++(\shift+0.8, 0);
	\draw (t) -- node[nix, below right]{$n_3n_4$} ++(0, -1);
	\node at (\posH +1.5*\Nshift+0.4, -2*\spacing-1.4) {(h)};

	\node[RT, label=above left:{$\vec W_1$}] (u1) at (\posI,-3*\spacing){};
	\node[RT, label=above left:{$\vec W_2$}] (u2) at (\posI+\Nshift,-3*\spacing){};
	\node[T, label=above:{$\vec x_2$}] (t) at (\posI+2*\Nshift,-3*\spacing){};
	\draw (u1) -- node[nix, below right]{$n_1$} ++(0,-1);
	\draw (u1) -- node[nix, below]{$r_1$} ++(\shift, 0);
	\draw (u2) -- node[nix, below right]{$n_2$} ++(0, -1);
	\draw (u2) -- node[nix, below]{$r_2$} ++(\shift, 0);
	\draw (t) -- node[nix, below right]{$n_3$} ++(0, -1);
	\draw (t) -- node[nix, below right]{$n_4$} ++(1, 0);
	\node at (\posI +\Nshift+0.5, -3*\spacing-1.4) {(i)};

	\node[GRT, label=above left:{$\vec W_1$}] (u1) at (\posJ,-3*\spacing){};
	\node[GRT, label=above left:{$\vec W_2$}] (u2) at (\posJ+\Nshift,-3*\spacing){};
	\node[T, label=above:{$\hat M_{(1,2)}(\vec x_2)$}] (t) at (\posJ+3*\Nshift,-3*\spacing){};
	\draw (u1) -- node[nix, below right]{$n_1$} ++(0,-1);
	\draw (u1) -- node[nix, below]{$r_1$} ++(\shift, 0);
	\draw (u2) -- node[nix, below right]{$n_2$} ++(0, -1);
	\draw (u2) -- node[nix, below]{$r_2$} ++(0.8*\shift, 0);
	\draw (t) -- node[nix, below]{$r_2n_3$} ++(-\shift, 0);
	\draw (t) -- node[nix, below right]{$n_4$} ++(1, 0);
	\node at (\posJ +1.5*\Nshift+0.5, -3*\spacing-1.4) {(j)};

	\node[GRT, label=above left:{$\vec W_1$}] (u1) at (\posK,-4*\spacing){};
	\node[GRT, label=above left:{$\vec W_2$}] (u2) at (\posK+\Nshift,-4*\spacing){};
	\node[RT, label=above left:{$\hat M_{(1,2)}(\vec W_3)$}] (u3) at (\posK+3*\Nshift,-4*\spacing){};
	\node[DT, label=above left:{$\vec S_2$}] (s) at (\posK+4*\Nshift,-4*\spacing){};
	\node[LT, label=above right:{$\vec V_2$}] (v) at (\posK+5*\Nshift,-4*\spacing){};
	\draw (u1) -- node[nix, below right]{$n_1$} ++(0,-1);
	\draw (u1) -- node[nix, below]{$r_1$} ++(\shift, 0);
	\draw (u2) -- node[nix, below right]{$n_2$} ++(0, -1);
	\draw (u2) -- node[nix, below]{$r_2$} ++(0.8*\shift, 0);
	\draw (u3) -- node[nix, below]{$r_2n_3$} ++(-\shift, 0);
	\draw (u3) -- node[nix, below]{$r_3$} ++(\shift, 0);
	\draw (s) -- node[nix, below]{$r_3$} ++(\shift, 0);
	\draw (v) -- node[nix, below right]{$n_4$} ++(0, -1);
	\node at (\posK +2.5*\Nshift, -4*\spacing-1.4) {(k)};

	\node[GRT, label=above left:{$\vec W_1$}] (u1) at (\posL,-5*\spacing){};
	\node[GRT, label=above left:{$\vec W_2$}] (u2) at (\posL+\Nshift,-5*\spacing){};
	\node[RT, label=above left:{$\hat M_{(1,2)}(\vec W_3)$}] (u3) at (\posL+3*\Nshift,-5*\spacing){};
	\node[T, label=above:{$\vec U_4$}] (u4) at (\posL+4*\Nshift,-5*\spacing){};
	\draw (u1) -- node[nix, below right]{$n_1$} ++(0,-1);
	\draw (u1) -- node[nix, below]{$r_1$} ++(\shift, 0);
	\draw (u2) -- node[nix, below right]{$n_2$} ++(0, -1);
	\draw (u2) -- node[nix, below]{$r_2$} ++(0.8*\shift, 0);
	\draw (u3) -- node[nix, below]{$r_2n_3$} ++(-\shift, 0);
	\draw (u3) -- node[nix, below]{$r_3$} ++(\shift, 0);
	\draw (u4) -- node[nix, below right]{$n_4$} ++(0, -1);
	\node at (\posL +2*\Nshift, -5*\spacing-1.4) {(l)};
	
	\node[RT, label=above left:{$\vec W_1$}] (u1) at (\posM,-5*\spacing){};
	\node[RT, label=above left:{$\vec W_2$}] (u2) at (\posM+\Nshift,-5*\spacing){};
	\node[RT, label=above left:{$\vec W_3$}] (u3) at (\posM+2*\Nshift,-5*\spacing){};
	\node[T, label=above:{$\vec W_4$}] (u4) at (\posM+3*\Nshift,-5*\spacing){};
	\draw (u1) -- node[nix, below right]{$n_1$} ++(0,-1);
	\draw (u1) -- node[nix, below]{$r_1$} ++(\shift, 0);
	\draw (u2) -- node[nix, below right]{$n_2$} ++(0, -1);
	\draw (u2) -- node[nix, below]{$r_2$} ++(\shift, 0);
	\draw (u3) -- node[nix, below right]{$n_3$} ++(0,-1);
	\draw (u3) -- node[nix, below]{$r_3$} ++(\shift, 0);
	\draw (u4) -- node[nix, below right]{$n_4$} ++(0, -1);
	\node at (\posM +1.5*\Nshift, -5*\spacing-1.4) {(m)};
\end{tikzpicture}
\caption[Step by step depiction of the TT-SVD by example]{Step by step depiction of the TT-SVD by example for an order $4$ tensor.}
\label{fig:TTCreation}
\end{figure}
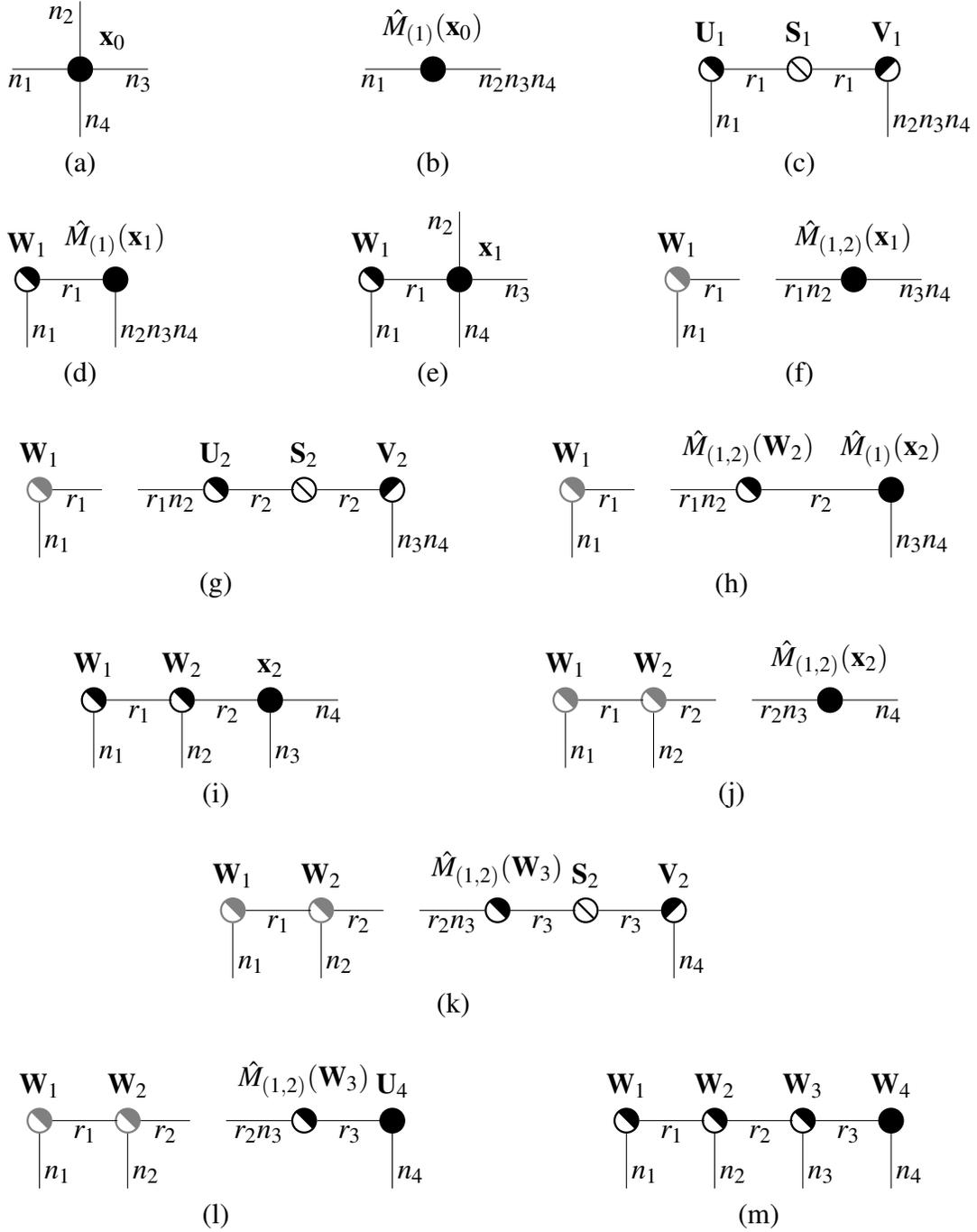

{\bf Tensor Train Singular Value Decomposition (TT-SVD)}\\
The intuition of the TT-SVD is that in every step a (matrix) SVD is performed to detach one open mode from the tensor. Figure \ref{fig:TTCreation} shows this process step by step for an order four tensor and is frequently referred to in the following description. The TT-SVD starts by calculating an SVD of the matricization of $\vec x = \vec x_0$, where all modes but the first one are combined (Fig. \ref{fig:TTCreation} (a)-(c))
\begin{align}
			\vec U_1 \vec S _1 \vec V_1^T &:= \text{SVD}\left( \hat M_{\{1\}} (\vec x_0)\right) \ ,
\end{align}
with $\vec U_1 \in \mathbb{R}^{n_1 \times r_1}$, $\vec S_1 \in \mathbb{R}^{r_1 \times r_1}$, $\vec V_1^T \in \mathbb{R}^{r_1 \times (n_2 \cdot  \ldots  \cdot n_d)}$. The dimension $r_1$ is equal to the rank of $\hat M_{\{1\}} \left( \vec x_0 \right)$. The resulting matrices $\left( \vec S_1 \vec V_1^T\right) $ and $\vec U_1$ are each dematricized, which is trivial for $\vec U_1$ in the first step (Fig. \ref{fig:TTCreation} (d)-(e)).
\begin{align}
    \vec W_1 &:= \hat M^{-1}(\vec U_1) & \vec W_1 \in  \mathbb{R}^{n_1 \times r_1}\\
    \vec x_1 &:= \hat M^{-1} \left( \vec S_1 \vec V_1^T \right) & \vec x_1 \in  \mathbb{R}^{r_1 \times n_2 \times  \ldots \times n_d}
\end{align}
Note that there holds
\begin{align}
	\vec W_1 \circ \vec x_1 = \hat M^{-1}\left( \vec U_1 \vec S_1 \vec V_1^T \right) = \vec x_0 = \vec x \ .
\end{align}
In the next step a matricization of the newly acquired tensor $\vec x_1$ is performed. The first dimension of the matricization is formed by the first two modes of $\vec x_1$, corresponding to the new dimension introduced by the prior SVD and the second original dimension. The second dimension of the matricization is formed by all remaining modes of $\vec x_1$ (Fig. \ref{fig:TTCreation} (f)). From this matricization another SVD is calculated (Fig. \ref{fig:TTCreation} (g))
\begin{align}
	\vec U_2 \vec S _2 \vec V_2^T &:= \text{SVD}\left( \hat M_{\{1,2\}} (\vec x_1)\right) \ ,
\end{align}
with $\vec U_2 \in \mathbb{R}^{(r_1 \cdot n_2) \times r_2}$, $\vec S_2 \in \mathbb{R}^{r_2 \times r_2}$, $\vec V_2^T \in \mathbb{R}^{r_2 \times (n_3 \cdot  \ldots  \cdot n_d)}$. As in the first step $\vec U_2$ and $\left( \vec S_2 \vec V_2^T\right)$ are then dematricized (Fig. \ref{fig:TTCreation} (i)),
\begin{align}
			\vec W_2 &:= \hat M^{-1}(\vec U_2) & \vec W_2 \in \mathbb{R}^{r_1 \times n_2 \times r_2}\\
			\vec x_2 &:= \hat M^{-1}(\vec S_2 \vec V_2^T) & \vec x_2 \in \mathbb{R}^{r_2 \times n_3 \times \ldots \times n_d} 
\end{align}
and again there holds
\begin{align}
			\vec W_2 \circ \vec x_2 = \hat M^{-1}\left( \vec U_2 \vec S_2 \vec V_2^T \right) &= \vec x_1\\
			\Rightarrow \vec W_1 \circ \vec W_2 \circ \vec x_2 &= \vec x \ .
\end{align}

The obtained rank $r_2$ is equal to the rank of the matricization $\hat M_{\{1,2\}}\left( \vec x\right)$, which can be shown as follows. First note that $\vec Q := \hat M_{\{1,2\}}\left( \vec U_1 \circ \vec U_2 \right) \in \mathbb{R}^{(n_1 \cdot n_2) \times r_2}$ is an orthogonal matrix, because
\begin{align}
			\left( \vec Q \vec Q^T\right) [j,j'] =& \sum_{i} \hat M_{\{1,2\}}\left( \vec W_1 \circ \vec W_2 \right)[i,j] \ \hat M_{(1,2)}\left( \vec W_1 \circ \vec W_2 \right)[i,j'] \\
			&= \sum_{i_1,i_2,k,k'} \vec W_1[i_1,k] \ \vec W_2[k,i_2,j] \ \vec W_1[i_1,k'] \ \vec W_2[k',i_2,j']\\
			&= \sum_{i_2,k,k'} \underbrace{\sum_{i_1} \vec W_1[i_1,k] \ \vec W_1[i_1,k']}_{\vec I[k,k']}  \ \vec W_2[k,i_2,j] \ \vec W_2[k',i_2,j']\\
			&= \sum_{i_2,k}\vec W_2[k,i_2,j] \ \vec W_2[k,i_2,j']\\
			&= \vec I [j, j'] \ .
\end{align}
Then consider that
\begin{align}
			 \hat M_{\{1,2\}}\left( \vec x \right) &= \hat M_{\{1,2\}}\left( \vec W_1 \circ \vec W_2 \circ \vec x_3 \right) \\
			 &= \hat M_{\{1,2\}}\left( \vec W_1 \circ \vec W_2 \right) \hat M_{\{1\}}\left(\vec x_2 \right) \\
			 &= \vec Q \vec S_2 \vec V_2^T\label{eq:ttSvdRankProof}
\end{align}
holds. This is a valid SVD of $\hat M_{\{1,2\}}( \vec x)$ and since the matrix of singular values is unique it follows that in fact $\text{rank}( \hat M_{\{1,2\}}( \vec x)) = r_2$.\\

This procedure is continued for a total of $d-2$ steps and in each step the order of $\vec x_i \in \mathbb{R}^{r_i \times n_{i+1} \times \ldots \times n_d}$ shrinks by one. Furthermore there holds
\begin{align}
	\vec W_1 \circ \vec W_2 \circ \ldots \circ \vec W_i \circ \vec x_i &= \vec x \,
\end{align}
$\vec W_i \in \mathbb{R}^{r_{i-1} \times n_i \times r_i}$ and $r_i = \text{rank}\left( \hat M_{\{1, \ldots, i\}}\left( \vec x\right)  \right)$ in every step. The $(d-1)$-st step (Fig. \ref{fig:TTCreation} (k)) 
\begin{align}
	\vec U_{d-1} \vec S_{d-1} \vec V_{d-1}^T &= \text{SVD}\left( \hat M_{\{1,2\}} (\vec x_{d-1})\right) \ ,
\end{align}
with $\vec U_{d-1} \in \mathbb{R}^{(r_{d-2} \cdot n_{d-1}) \times r_{d-1}}$, $\vec S_{d-1} \in \mathbb{R}^{r_{d-1}\times r_{d-1}}, \vec V_{d-1}^T \in  \mathbb{R}^{r_{d-1} \times n_d}$, is special since the de-matricization of $(\vec S _{d-1} \vec V_{d-1}^T)$, yields an order two tensor that is named $\vec W_d$ instead of $\vec x_{d}$ (Fig. \ref{fig:TTCreation} (l)-(m))
\begin{align}
	\vec W_{d-1} &= \hat M^{-1} (\vec U_{d-1}) & \vec W_{d-1} \in \mathbb{R}^{r_{d-2} \times n_{d-1} \times r_{d-1}}\\
	\vec W_{d} &= \vec S_{d-1} \vec V_{d-1}^T & \vec x_{d} \in \mathbb{R}^{r_{d-1} \times n_d} \ .
\end{align} 
Finally
\begin{align}
	\vec W_1 \circ \vec W_2 \circ \ldots \circ \vec W_{d-1} \circ \vec W_d &= \vec x
\end{align}
is a valid TT representation of $\vec x$ with TT-rank $\vec r = (r_1, \ldots, r_{d-1})$, whose entries $r_i = \operatorname{rank} (\hat M_{\{1, \ldots, i\}}\left( \vec x\right))$ are exactly the ranks of the matricizations as asserted.\\

The same algorithm can also be used to calculate a rank $\vec r^* = (r_1^*, \ldots, r_{d-1}^*)$ approximation of a tensor $\vec x \in \mathbb{R}^{n_1 \times \ldots \times n_d}$ with TT-rank $\vec r \succeq \vec r^*$. To this end the normal SVDs are replaced by truncated rank $r_i^*$ SVDs, yielding a tensor $\vec x^*$ of TT-rank $\vec r^*$. In contrast to the matrix case, $\vec x^*$ is in general not the best rank $\vec r^*$ approximation of $\vec x$. However as shown by \cite{oseledets2011tensor}, it is a quasi best approximation with
\begin{align}
	 \lVert \vec x-\vec x^* \rVert_2 \leq \sqrt{d-1} \min_{\vec y \ :\ \operatorname{TT-rank}(\vec y) \preceq \vec r^*} \left(\lVert \vec x - \vec y \rVert \right) \ .
\end{align}
The computational complexity of the TT-SVD is dominated by the $d-1$ matrix singular value decompositions, with all other contributions being asymptotically negligible. With $n := \max(n_1, \ldots, n_d)$ and $r := \max(r_1, \ldots r_d)$ the cost scales as $\mathcal{O}(n^{d+1} + \sum_{i=1}^{d-1} r^2 n^{d-i}) \subset \mathcal{O}(dn^{d+1})$, i.e. still exponential in the order. This is somewhat expected because there are in general $n^d$ entries in the original tensor that have to be considered. Unfortunately $\vec x$ being sparse or otherwise structured incurs no dramatic change because the structure is generally lost after the first SVD.\\

Apart from the generalized singular value decomposition the TT format offers several further beneficial properties. In particular it is able to break the curse of dimensionality, in the sense that the  storage complexity of a tensor $\vec x \in \mathbb{R}^{n_1 \times \ldots \times n_d}$ with TT-rank $\vec r =(r_1, \ldots, r_{d-1})$ in a minimal TT representation scales as $\mathcal{O}(dnr^2)$, i.e.\ linearly in the order. Here $n := \max_i(n_1, \ldots, n_d)$ and $r := \max_i(r_1, \ldots, r_{d-1})$. Additionally also the computational complexity of common operations as additions and scalar products, scale only linearly in the order for fixed ranks, see \cite{oseledets2011tensor, hackbusch2012tensorBuch}. Another desirable property is that the set of tensors with rank at most $\vec r$ form a closed set and as shown by \citet{holtz2012manifolds} the set of tensor with exact rank $\vec r$ forms a smooth manifold, allowing the application of Riemannian optimization techniques \cite{kressner2016preconditioned, steinlechner2016riemannian} and dynamical low rank approximation \cite{lubich2013dynamical, lubich2015time}, see also the review article \cite{bachmayr2016tensor}. Especially for numerical applications these properties made the tensor train one of, if not the, most popular tensor decomposition of recent years.

\section{Randomized SVD for higher order tensors}
As shown in the previous section calculating a low rank representation or approximation of a given higher order tensor is a challenging task, as the complexity of the tensor train SVD (TT-SVD) scales exponentially in the order. For dense tensors this is of course somewhat expected as there is an exponential number of entries that have to be incorporated. Nevertheless also for sparse and structured matrices the two decomposition techniques exhibit an exponential scaling. In this section we look at randomized methods for the calculation of approximate matrix factorizations. For sparse or structured matrices these techniques allow for a very efficient calculation of common matrix factorizations such as the SVD or QR decomposition, while offering rigorous stochastic error bounds. In the second part of this section we apply these results to formulate randomized TT-SVD algorithms. We show that there hold stochastic error bounds similar to the matrix setting. We also show that this randomized TT-SVD has only linear complexity with respect to the order when applied to sparse tensors.

\subsection{Randomized SVD for matrices}
\label{sec:randMatrixDecomp}

Randomized techniques for the calculation of SVD or QR factorizations of matrices have been proposed many times in the literature. However it was only recently that, thanks to the application of new results from {\em random matrix theory} these procedures could be analyzed rigorously. We start this section by presenting some results from the work of \citet{halko2011finding}, which will provide a solid basis for the randomized tensor factorization methods of the second part of this section. In this part we restrict ourself to standard Gaussian random matrices, i.e. matrices whose entries are i.i.d.\ standard Gaussian random variables. The usage of structured random matrices is discussed in section \ref{sec:conclusions}. \\ 

In the formulation of \citet{halko2011finding}, the basis of all decompositions is a randomized method to calculate an approximate low rank subspace projection
\begin{align} \label{eq:rangeProjection}
\vec A \approx \vec Q \vec Q^T \vec A
\end{align}
where $\vec A \in \mathbb{R}^{n \times m}$ is a given matrix and $\vec Q \in \mathbb{R}^{n \times s}$ is an orthogonal matrix approximately spanning the range of $\vec A$. Here $s = r + p$, where $r$ is the desired rank and $p$ is an oversampling parameter. With this projection at hand numerous different low rank decompositions can be calculated deterministically at low costs. For example the singular value decomposition of $\vec A$ can be calculated by forming $\vec B := \vec Q^T \vec A$ and calculating the deterministic SVD $\vec U \vec S \vec V^T = \vec B$ of $\vec B$. Using $\tilde{\vec U} = \vec Q \vec U$
\begin{align}
\tilde{\vec U} \vec S \vec V^T = \vec Q \vec Q^T \vec A \approx \vec A
\end{align}
is an approximate SVD of $\vec A$ containing only the approximation error incurred by the subspace-projection. The computational costs of the involved operations scale as $\mathcal{O}\left( s T_{\text{mult}} + s^2(m + n) \right)$, where $T_{\text{mult}}$ is the cost to calculate the matrix-vector product with $\vec A$, which is $\mathcal{O}(mn)$ for a general matrix but can be much lower for structured or sparse matrices. In a similar way other matrix factorizations can also be computed with low costs if the projection \eqref{eq:rangeProjection} is given.\\

The main challenge is the calculation of the approximate range $\vec Q$ through random techniques. For this \cite{halko2011finding} present the following prototype algorithm. Given a matrix $\mathbf{A} \in \mathbb{R}^{n_1 \times n_2}$.
\begin{lstlisting}[caption=Randomized range approximation, mathescape=true, label=list:randSVD]
Input: $\vec A,\  r,\  p$ Output: $\vec Q$
Create a standard Gaussian random matrix $\mathbf{G} \in \mathbb{R}^{n_2 \times (r+p)}$
Calculate the intermediate matrix $\vec B := \mathbf{A} \mathbf{G} \in \mathbb{R}^{n_1 \times s}$.
Compute the factorization $\vec Q \vec R = \vec B$.
\end{lstlisting}

The following theorem proves that the $\mathbf{Q}$ obtained in this manner is indeed an approximation of the range of $\mathbf{A}$ in the sense of \eqref{eq:rangeProjection}.
\begin{theorem}[\citet{halko2011finding}] \label{theo:randMatrixError}
Given $\vec A \in \mathbb{R}^{m \times n}$ and $s = r + p$ with $p \ge 2$. For the projection $\vec Q$ obtained by procedure \ref{list:randSVD} there holds the following error bounds.
\begin{equation}
\| \mathbf{A} - \mathbf{Q} \mathbf{Q}^T \mathbf{A}  \| \leq 
\left[1+11 \sqrt{(r + p) \cdot \min(m,n)}\right] \sigma_{r+1}
\end{equation}
with probability at least $ 1-6 p^{-p} $ and for $p\ge 4$ and any $u,t \ge 1$
\begin{equation}
\| \mathbf{A} - \mathbf{Q} \mathbf{Q}^T \mathbf{A}  \| \leq 
\left[1+t \sqrt{\frac{12r}{p}}\right] \left(\sum_{k>r} \sigma_k^2\right)^{1/2} + u t \frac{e\sqrt{r+p}}{p+1}\sigma_{r+1}
\end{equation}
with probability at least $1 - 5t^{-p}-2e^{-u^2/2}$.
\end{theorem}

Let us highlight furthermore that for the operator norm, we have that $ \sigma_{r+1} = \mbox{ inf }_{\mbox{rank } \mathbf{B} = r }  \| \mathbf{A} - \mathbf{B} \|_{op} $ and for the Frobenius norm there holds 
$$ 
\left( \sum_{k > r } \sigma_k^2 \right)^{\frac{1}{2}} = \inf_{\text{rank}(\mathbf{B}) \leq r}  \| \mathbf{A} - \mathbf{B} \| \ .
$$

\subsection{Randomized TT-SVD}
\label{sec:randTensorDecomp}

In this section we show how the same idea of the randomized range approximation for matrices can be used to formulate a randomized algorithm that calculates an approximate TT-SVD of arbitrary tensors. We show that stochastic error bounds analogous to the matrix case can be obtained. Furthermore we show that for sparse tensors this randomized TT-SVD can be calculated in \emph{linear} complexity with respect to the order of the tensor, instead of the exponential complexity of the deterministic TT-SVD.\\

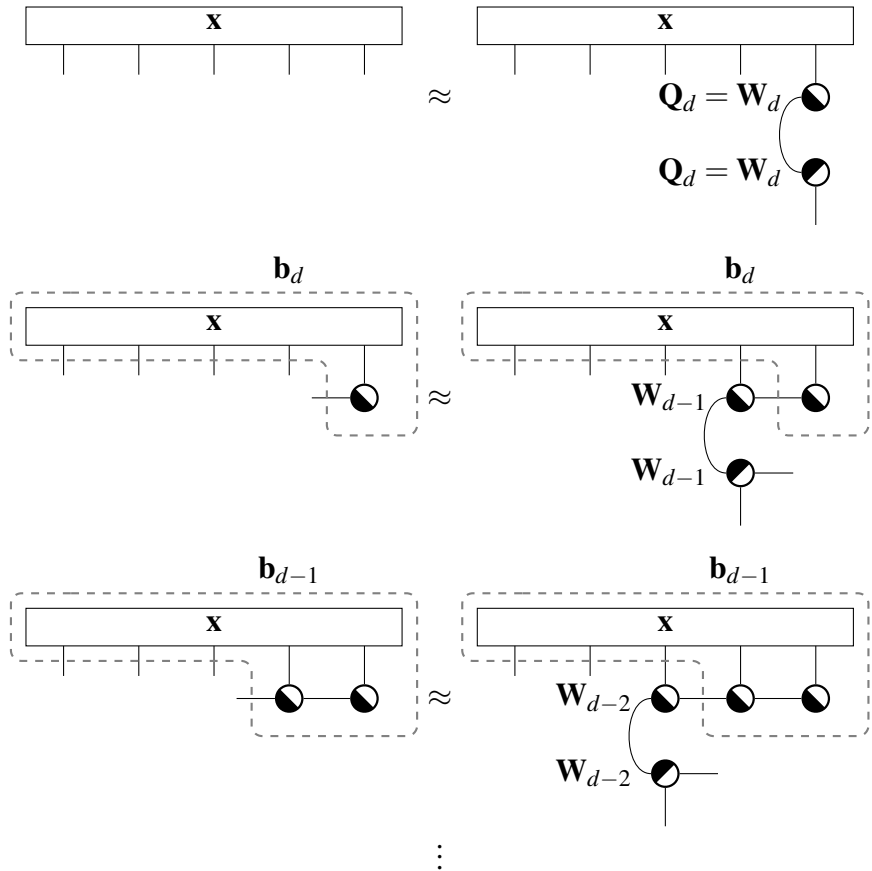
\begin{figure}
\centering	
\begin{tikzpicture}
\begin{scope}[shift={(0,0)}]
\draw[dashed, thick, white, rounded corners, loop]
		(0.1,1.4) -- (4.7,1.4) -- (4.7,-.5) -- (3.5,-.5) -- (3.5,.5) -- (-.7,.5) -- (-.7,1.4) -- (0.1,1.4);
	\draw (0,.3) -- (0,.7);
	\draw (1,.3) -- (1,.7);
	\draw (2,.3) -- (2,.7);
	\draw (3,.3) -- (3,.7);
	\draw (4,.3) -- (4,.7);
	\draw (-.5,1.2) rectangle (4.5,.7);
	\node[nix] (A) at (2,1){$\vec x$};
    
	\node[nix] at (5,0){$\approx$};
    
\begin{scope}[shift={(6,0)}]
\draw[dashed, thick, white, rounded corners, loop]
		(0.1,1.4) -- (4.7,1.4) -- (4.7,-.5) -- (3.5,-.5) -- (3.5,.5) -- (-.7,.5) -- (-.7,1.4) -- (0.1,1.4);
	\node[LLT, label=above right:{$\vec Q_d = \vec W_d \ $}] (u5t) at (4,0){};
	\node[LT, label=above left:{$\vec Q_d = \vec W_d \ $}] (u5) at (4,-1){} edge[bend left=90] (u5t);
		\draw (0,.3) -- (0,.7);
		\draw (1,.3) -- (1,.7);
		\draw (2,.3) -- (2,.7);
		\draw (3,.3) -- (3,.7);
	\draw (u5t) -- ++(0,.7);
	\draw (u5) -- ++(0,-.7);
	\draw (-.5,1.2) rectangle (4.5,.7);
	\node[nix] (A) at (2,1){$\vec x$};
	\node[nix] (B) at (3,1.7){$\textcolor{white}{\vec B_d}$};
\end{scope}
\end{scope}

\begin{scope}[shift={(0,-4)}]
	\node[LLT] (u5t) at (4,0){};
		\draw (0,.3) -- (0,.7);
		\draw (1,.3) -- (1,.7);
		\draw (2,.3) -- (2,.7);
		\draw (3,.3) -- (3,.7);
	\draw (u5t) -- ++(0,.7);
	\draw (u5t) -- ++(-.7,0);
	\draw (-.5,1.2) rectangle (4.5,.7);
	\node[nix] (A) at (2,1){$\vec x$};
	\draw[dashed, thick, gray, rounded corners, loop]
		(0.1,1.4) -- (4.7,1.4) -- (4.7,-.5) -- (3.5,-.5) -- (3.5,.5) -- (-.7,.5) -- (-.7,1.4) -- (0.1,1.4);
	\node[nix] (B) at (3,1.7){$\vec b_d$};
    
	\node[nix] at (5,0){$\approx$};
    
    \begin{scope}[shift={(6,0)}]
		\node[LLT, label=above right:{$\vec W_{d-1}\ $}] (u4t) at (3,0){};
		\node[LT, label=above left:{$\vec W_{d-1}\ $}] (u4) at (3,-1){} edge[bend left=90] (u4t);
		\node[LLT] (u5t) at (4,0){} edge (u4t);
			\draw (0,.3) -- (0,.7);
			\draw (1,.3) -- (1,.7);
			\draw (2,.3) -- (2,.7);
		\draw (u5t) -- ++(0,.7);
      \draw (u4t) -- ++(0,.7);
		\draw (u4) -- ++(0,-.7);
		\draw (u4) -- ++(.7,0);
		\draw (-.5,1.2) rectangle (4.5,.7);
		\node[nix] (A) at (2,1){$\vec x$};
		\draw[dashed, thick, gray, rounded corners, loop]
			(0.1,1.4) -- (4.7,1.4) -- (4.7,-.5) -- (3.5,-.5) -- (3.5,.5) -- (-.7,.5) -- (-.7,1.4) -- (0.1,1.4);
		\node[nix] (B) at (3,1.7){$\vec b_d$};
	\end{scope}
\end{scope}

\begin{scope}[shift={(0,-8)}]
	\node[LLT] (u4t) at (3,0){};
	\node[LLT] (u5t) at (4,0){} edge (u4t);
		\draw (0,.3) -- (0,.7);
		\draw (1,.3) -- (1,.7);
		\draw (2,.3) -- (2,.7);
	\draw (u5t) -- ++(0,.7);
	\draw (u4t) -- ++(0,.7);
	\draw (u4t) -- ++(-.7,0);
	\draw (-.5,1.2) rectangle (4.5,.7);
	\node[nix] (A) at (2,1){$\vec x$};
	\draw[dashed, thick, gray, rounded corners, loop]
		(0.1,1.4) -- (4.7,1.4) -- (4.7,-.5) -- (2.5,-.5) -- (2.5,.5) -- (-.7,.5) -- (-.7,1.4) -- (0.1,1.4);
	\node[nix] (B) at (3,1.7){$\vec b_{d-1}$};
    
    \node[nix] at (5,0){$\approx$};
    
\begin{scope}[shift={(6,0)}]
	\node[LLT, label=above right:{$\vec W_{d-2}\ $}] (u3t) at (2,0){};
	\node[LT, label=above left:{$\vec W_{d-2}\ $}] (u3) at (2,-1){} edge[bend left=90] (u3t);
	\node[LLT] (u4t) at (3,0){} edge (u3t);
	\node[LLT] (u5t) at (4,0){} edge (u4t);
		\draw (0,.3) -- (0,.7);
		\draw (1,.3) -- (1,.7);
	\draw (u5t) -- ++(0,.7);
	\draw (u4t) -- ++(0,.7);
	\draw (u3t) -- ++(0,.7);
	\draw (u3) -- ++(0,-.7);
	\draw (u3) -- ++(.7,0);
	\draw (-.5,1.2) rectangle (4.5,.7);
	\node[nix] (A) at (2,1){$\vec x$};
	\draw[dashed, thick, gray, rounded corners, loop]
			(0.1,1.4) -- (4.7,1.4) -- (4.7,-.5) -- (2.5,-.5) -- (2.5,.5) -- (-.7,.5) -- (-.7,1.4) -- (0.1,1.4);
		\node[nix] (B) at (3,1.7){$\vec b_{d-1}$};
\end{scope}
\end{scope}	

\node[nix] at (5,-10){$\vdots$};

\end{tikzpicture}
\caption{Iterative construction of the tensor $\vec b_i$ by subsequent range approximations.}
\label{fig:RandTTSVD}
\end{figure}

The idea of our randomized TT-SVD procedure is to calculate nested range approximations increasing by one mode at a time. The corresponding projector is composed of separated orthogonal parts, which are calculated using procedure \ref{list:randSVD}. This is visualized in figure \ref{fig:RandTTSVD}. These orthogonal parts will become the component tensors $\vec W_2, \ldots, \vec W_d$ of the final TT decomposition. The first component tensor $\vec W_1$ is given by contracting the initial $\vec x$ with all orthogonal components, i.e.\ 
$$
\vec W_1 = \vec x \circ_{(2, \ldots, d), (2, \ldots, d)} \left( \vec W_2 \circ \ldots \circ \vec W_d \right)
$$
The exact procedure calculating the orthogonal components and this final contraction is given in listing \ref{list:randTTSVD}.

\begin{lstlisting}[caption=Randomized TT-SVD, label=list:randTTSVD, mathescape=true]
Input: $\vec x$, Output: $\vec W_1, \ldots, \vec W_d$

Set $\vec b_{d+1} := \vec x$
For $j=d, \ldots, 2$:
  Create a Gaussian random tensor $\vec g \in \mathbb{R}^{s_{j-1} \times n_1 \times \ldots \times n_{j-1}}$ 
  Calculate $\vec a_{j} := \vec g \circ_{(2, \ldots j), (1, \ldots, j-1)} \vec b_{j+1}$ 
  Calculate the factorization $\vec R_{j} \vec Q_{j} := \hat M_{\{1\}} \left( \vec a_{j} \right)$ 
  Set $\vec W_{j} := \hat M^{-1} \left( \vec Q_{j} \right)$ 
  if j = d:
    Calculate $\vec b_{j} = \vec b_{j+1} \circ_{(j), (2)} \vec W_{j}$
  else
    Calculate $\vec b_{j} = \vec b_{j+1} \circ_{(j, j+1), (2, 3)} \vec W_{j}$
    
Set $\vec W_1 = \vec b_2$
\end{lstlisting}

\begin{figure}
\centering	
\begin{tikzpicture}
\begin{scope}[shift={(0,0)}]
	\node[LLT] (u2t) at (1,0){};
	\node[LLT] (u3t) at (2,0){} edge (u2t);
	\node[LLT] (u4t) at (3,0){} edge (u3t);
	\node[LLT] (u5t) at (4,0){} edge (u4t);
	\node[LT] (u2) at (1,-1){} edge[bend left=90] (u2t);
	\node[LT] (u3) at (2,-1){} edge (u2);
	\node[LT] (u4) at (3,-1){} edge (u3);
	\node[LT] (u5) at (4,-1){} edge (u4);
		\draw (0,-1.7) -- (0,.7);
	\draw (u5t) -- ++(0,.7);
	\draw (u4t) -- ++(0,.7);
	\draw (u3t) -- ++(0,.7);
	\draw (u2t) -- ++(0,.7);
	\draw (u2) -- ++(0,-.7);
	\draw (u3) -- ++(0,-.7);
	\draw (u4) -- ++(0,-.7);
	\draw (u5) -- ++(0,-.7);
	\draw (u3) -- ++(.7,0);
	\draw (-.5,1.2) rectangle (4.5,.7);
	\node[nix] (A) at (2,1){$\vec x$};
	\draw[dashed, thick, gray, rounded corners, loop]
			(0.1,1.4) -- (4.7,1.4) -- (4.7,-.5) -- (-.7,-.5) -- (-.7,1.4) -- (0.1,1.4);
		\node[nix] (B) at (3,1.7){$\vec B_2$};
\end{scope}	
\node[nix] at (5,0){$=$};
\begin{scope}[shift={(6,0)}]
	\node[LLT] (u2t) at (1,0){};
	\node[LLT] (u3t) at (2,0){} edge (u2t);
	\node[LLT] (u4t) at (3,0){} edge (u3t);
	\node[LLT] (u5t) at (4,0){} edge (u4t);
	\node[LT] (u2) at (1,-1){} edge[bend left=90] (u2t);
	\node[LT] (u3) at (2,-1){} edge (u2);
	\node[LT] (u4) at (3,-1){} edge (u3);
	\node[LT] (u5) at (4,-1){} edge (u4);
		\draw (0,-1.7) -- (0,.7);
	\draw (u5t) -- ++(0,.7);
	\draw (u4t) -- ++(0,.7);
	\draw (u3t) -- ++(0,.7);
	\draw (u2t) -- ++(0,.7);
	\draw (u2) -- ++(0,-.7);
	\draw (u3) -- ++(0,-.7);
	\draw (u4) -- ++(0,-.7);
	\draw (u5) -- ++(0,-.7);
	\draw (u3) -- ++(.7,0);
	\draw (-.5,1.2) rectangle (4.5,.7);
	\node[nix] (A) at (2,1){$\vec x$};
	\draw[dashed, thick, gray, rounded corners, loop]
			(0.1,.5) -- (4.5,.5) -- (4.5,-1.5) -- (-.5,-1.5) -- (-.5,.5) -- (0.1,.5);
		\node[nix] (B) at (5.1,-.5){$\hat P_{2, \ldots, d}$};
\end{scope}

\end{tikzpicture}
\caption{Depiction of the randomized TT-SVD as the action of the projection operator $\hat P_{2,\ldots,d}$}
\label{fig:projector}
\end{figure}
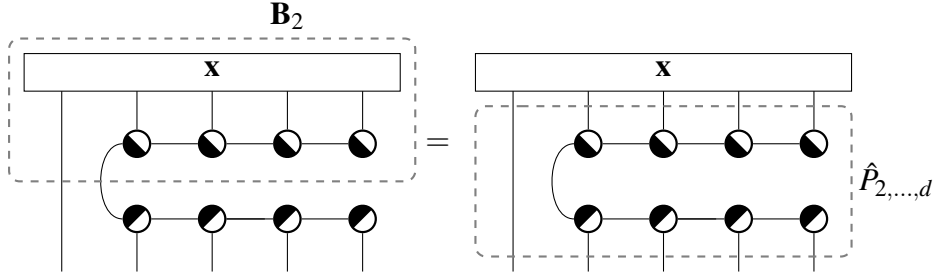

At the end of the procedure in listing \ref{list:randTTSVD}
$$ \vec x \approx \vec W_1 \circ \vec W_2 \circ \vec W_3 \circ \ldots \circ \vec W_d$$
is an approximate TT decomposition of rank $\vec s = (s_1, \ldots, s_{d-1})$. This final composition can also be given in terms of contractions with the orthogonal parts, i.e.\ 
\begin{align}
\vec x &\approx \vec W_1 \circ \vec W_2 \circ \vec W_3 \circ \ldots \circ \vec W_d\\
&= \left( \vec x \circ_{(2, \ldots, d), (2, \ldots, d)} \left( \vec W_2 \circ \ldots \circ \vec W_d \right) \right) \circ \vec W_2 \circ \vec W_3 \circ \ldots \circ \vec W_d\\
&= \vec x \circ_{(2, \ldots, d), (2, \ldots, d)} \left( \left( \vec W_2 \circ \ldots \circ \vec W_d \right) \circ_{1, 1} \left( \vec W_2 \circ \ldots \circ \vec W_d \right) \right) \\
&=: \hat P_{2, \ldots, d} ( \vec x)
\end{align}
where the effect of the orthogonal parts can also be seen as the action of an projector $\hat P_{2, \ldots, d}$. Note that since all parts are orthogonal this is indeed an orthogonal projector. This relation is visualized in figure \ref{fig:projector}. In the following it will be useful to also define the orthogonal projections
\begin{equation}
\hat P_{i, \ldots, d} (\vec x) := \vec x \circ_{(i, \ldots, d), (i, \ldots, d)} \left( \left( \vec W_i \circ \ldots \circ \vec W_d \right) \circ_{1, 1} \left( \vec W_i \circ \ldots \circ \vec W_d \right) \right) \ .
\end{equation}

The following theorem shows that there exists an stochastic error bound for this randomized TT-SVD.
\begin{theorem}[Error bound] \label{theo:RandTTSVDError}
	Given $\vec x \in \mathbb{R}^{n_1 \times \ldots \times n_d}$, $s= r + p$ with $p \geq 4$. For every $u,t \geq 1$ the error of the randomized TT-SVD, as given in listing \ref{list:randTTSVD}, fulfills
	\begin{equation} \label{eq:stochErrorBoundRndTTSvD}
		\|\vec x - P_{2..d}(\vec x)\| \le \sqrt{d-1}\; \eta(r, p) \min_{\text{TT-rank}(\vec y) \leq \vec r } \|\vec x - \vec y\|
	\end{equation}
	with probability at least $(1 - 5t^{-p}-2e^{-u^2/2})^{d-1}$. The parameter $\eta$ is given as
    \begin{equation}
    \eta = 1+t\sqrt{\frac{12r}{p}} + ut \frac{e\sqrt{r+p}}{p+1} \ .
    \end{equation}
     
\begin{proof}
For syntactical convenience let use define $\vec B_i := \left( \hat M_{\{1, \ldots, i-1\}} ( \vec b_i ) \right)^T$. Then as $\hat P_{2,\ldots,d}$ is an orthogonal projector we have
\begin{align} \label{eq:baseError}
\| \vec x -P_{2..d}(\vec x )\|^2 &= \|\vec x \|^2 - \|P_{2..d}(\vec x )\|^2 \\
&= \| \vec x \|^2 - \big<\vec B_2, \vec B_2\big>\\
&= \| \vec x \|^2 - \big<\vec Q_2 \vec B_3,  \vec Q_2 \vec B_3\big>\\
&= \| \vec x \|^2 - \big<\vec B_3, \vec Q_2^T \vec Q_2 \vec B_3\big> \ .
\end{align}
For all $2 \leq i \leq d$ there holds
\begin{align}
\big<\vec B_{i+1}, \vec Q_i^T \vec Q_i \vec B_{i+1}\big> &= \left<\vec B_{i+1}, \vec B_{i+1} - (\vec I - \vec Q_i^T \vec Q_i)\vec B_{i+1}\right> \\
&= \|\vec B_{i+1}\|^2 - \left<\vec B_{i+1}, (\vec I - \vec Q_i^T \vec Q_i)\vec B_{i+1}\right> \\
&= \|\vec B_{i+1}\|^2 - \|(\vec I - \vec Q_i^T \vec Q_i)\vec B_{i+1}\|^2 \\
&= \big<\vec B_{i+2}, \vec Q_{i+1}^T \vec Q_{i+1} \vec B_{i+2}\big> - \|(\vec I - \vec Q_i^T \vec Q_i)\vec B_{i+1}\|^2 \ ,
\end{align}
where we used that $\vec Q_i^T \vec Q_i$ is an orthogonal projector as well. Inserting this iteratively into \eqref{eq:baseError} gives
\begin{align}
\| \vec x - P_{2..d}(\vec x )\|^2 &= \|\vec x \|^2 - \|\vec B_{d+1}\|^2 + \sum_{i=2}^{d} \|(\vec I-\vec Q_i^T \vec Q_i)\vec B_{i+1}\|^2 \\
&= \sum_{i=2}^{d} \|(\vec I-\vec Q_i^T \vec Q_i)\vec B_{i+1}\|^2 \ ,
\end{align}
where we used that $\vec B_{d+1}$ as a matricization of $\vec b_{d+1} := \vec x$ has the same norm as $\vec x$ itself. As $\vec Q_i$ is obtained in the exact setting of theorem \ref{theo:randMatrixError}, we know that for all $i$
\begin{align}
& \ \ \ \ \|( \vec I - \vec Q_i \vec Q_i^T)\vec B_{i+1}\|^2 \\
&\leq \left[ \left(1+t \sqrt{\frac{12r}{p}}\right) \left(\sum_{k>r} \sigma_k^2(\vec B_{i+1})\right)^{1/2} + u t \frac{e\sqrt{r+p}}{p+1}\sigma_{r+1}(\vec B_{i+1}) \right]^2 \\ 
&\leq \left[ \left(1+t \sqrt{\frac{12r}{p}} + u t \frac{e\sqrt{r+p}}{p+1} \right) \left(\sum_{k>r} \sigma_k^2(\vec B_{i+1})\right)^{1/2} \right]^2 \\
&\leq \eta^2 \sum_{k>r} \sigma_k^2(\vec B_{i+1})
\end{align}
holds with probability at least $1 - 5t^{-p}-2e^{-u^2/2}$. Note that the singular values of $\vec B_{i+1}$ are the same as of $\hat M_{\{1, \ldots, i-1\}} (\hat P_{i+1..d}(\vec x))$, see e.g.\ figure \ref{fig:projector}. As shown by \citet{Hochstenbach20101053} the application of an orthogonal projection can only decrease the singular values. Thereby it follows that 
\begin{align}
\|( \vec I - \vec Q_i \vec Q_i^T)\vec B_{i+1}\|^2 &\leq \eta^2 \sum_{k>r} \sigma_k^2(\vec B_{i+1}) \\
	&= \eta^2\sum_{k>r} \sigma_k^2(\hat M_{\{1, \ldots, i-1\}} (\hat P_{i+1..d}(\vec x))) \\
	&\le \eta^2\sum_{k>r} \sigma_k^2(\hat M_{\{1, \ldots, i-1\}} (\vec x)) \\
	&\le \eta^2 \min_{\text{rank}(\hat M_{\{1, \ldots, i-1\}} (\vec y)) \leq r}\|\vec x - \vec y\|^2\\
    &\leq \eta^2 \min_{\text{TT-rank}(\vec y) \leq \vec r}\|\vec x - \vec y\|^2 \ .
\end{align}
As the random tensors $\vec g$ are sampled independently in each step, the combined probability that the above holds for all $i$ is at least $\rho \ge (1 - 5t^{-p}-2e^{-u^2/2})^{d-1}$, as asserted.
\end{proof}
\end{theorem}

Note that if the tensor $\vec x$ actually has TT-rank $\vec r$ or smaller, that is if \ $\min_{\text{TT-rank}(\vec y) \leq \vec r}\|\vec x - \vec y\| = 0$, then the randomized TT-SVD is exact with probability one. This follows directly from theorem \ref{theo:RandTTSVDError} by using $t \rightarrow \infty, u \rightarrow \infty$.\\

Using standard Gaussian random tensors the computational complexity of the randomized TT-SVD is bounded by $\mathcal{O}(d s n^d)$, which is very similar to the deterministic TT-SVD presented in section \ref{sec:ttFormat}. However, as we show in the following proposition \ref{prop:sparseComplexity}, for sparse tensors the complexity scales only \emph{linearly} in the order, which is a dramatic reduction compared to the exponential scaling of the deterministic TT-SVD.

\begin{proposition}
\label{prop:sparseComplexity}
Assume that $\vec x \in \mathbb{R}^{n_1 \times \ldots \times n_d}$ contains at most $N$ non-zero entries. Then the computational complexity of the randomized TT-SVD given in listing \ref{list:randTTSVD} scales as $\mathcal{O}(d(s^2 N +s^3 n))$.

\begin{proof}
First note that if $\vec x$ has at most $N$ non-zero entries then each $\vec b_i$ has at most $s_{i-1} N$ non-zero entries. The fact that $\vec x$ has at most $N$ non-zero entries implies that, independent of $j$, there are at most $N$ tuples $(k_1, \ldots, k_j)$ such that the sub-tensor $\vec x[k_1, \ldots, k_j, \cdot, \ldots, \cdot]$ is not completely zero. Now each $\vec b_i$ can be given as $\vec b_i = \vec x \circ_{(i, \ldots, d), (2, \ldots, d-i+1)} \left( \vec W_i \circ \ldots \circ \vec W_d \right)$. As any contraction involving a zero tensors results in a zero tensor, $\vec b_i$ is non zero only if the first $d-i+1$ modes take values according to the at most $N$ tuples. As there is only one further mode of dimension $s_{i-1}$ there can in in total be only $s_{i-1} N$ non-zero entries in $\vec b_i$. \\

Creating only the, at most $s_{j-1} s_j N$, entries of $\vec g$ actually needed to perform the product $\vec a_{j} := \vec g \circ_{(2, \ldots j), (1, \ldots, j-1)} \vec b_{j+1}$ this calculation can be done in $\mathcal{O}(s_{j-1} s_j N)$. Calculating the $QR$ of $\hat M_{(1)} \left( \vec a_j \right)$ has complexity $\mathcal{O}(s_{j-1}^2 n_j s_j)$. The involved (de-)matrification actually do not incur any computational costs. Finally the product $\vec b_{j} = \vec b_{j+1} \circ_{(j, j+1), (2, 3)} \vec W_{d-j}$ has complexity $\mathcal{O}(s_{j-1} s_j N)$. These steps have to be repeated $d-1$ times. Adding it all up this gives an asymptotic cost bounded by $\mathcal{O}(d(s^2 N + s^3 n))$, where $s := \max(s_1, \ldots, s_d)$.
\end{proof}
\end{proposition}

\section{Relation to the Alternating Least Squares (ALS) Algorithm}
\label{sec:als}

There is an interesting connection between the proposed randomized TT-SVD and the popular alternating least squares (ALS) algorithm, which is examined in this section. Most of this section is still work in progress but we consider sharing the ideas worthwhile nevertheless. The ALS itself is a general optimization algorithm, highly related to the very successful DMRG algorithm known in quantum physics. We provide only a minimal introduction and refer to the literature for an exhaustive treatment, see e.g.\ \cite{holtz2012alternating, espig2015convergence}. \\

The ALS is used to solve optimization problems on the set of tensors with fixed TT-rank $\vec r$, for general objective functionals $\mathcal{J} : \mathbb{R}^{n_1 \times \ldots \times n_d} \to \mathbb{R}$. The special case interesting in this work is
$$
\mathcal{J} (\mathbf{x}) := \| \mathbf{f} - \mathbf{x} \|^2  
$$
for a given tensor $\mathbf{f} \in \mathbb{R}^{n_1 \times \ldots \times n_d}$. The global optimum is then exactly the best rank $\vec r$ approximation
\begin{equation}
\mathbf{x}^* := \operatorname{argmin}_{\text{TT-rank}(\vec y) = \vec r} \| \vec y - \vec f \| \ . 
\end{equation}
Observing that the parametrization $\vec x = \tau \left( \vec W_1, \ldots, \vec W_d \right) = \vec W_1 \circ \ldots \circ \vec W_d$ is multi-linear in the parameters $\mathbf{W}_i \in \mathbb{R}^{r_{i-1} \times n_i \times r_i}$. Hence fixing all components $\mathbf{W}_1 $ except the $i$-th component $\mathbf{U}_i$, provides a parametrization of $\vec x$ which is linear in $\mathbf{U}_i \in \mathbb{R}^{r_{i-1} \times n_i \times r_i}$,
\begin{equation}
\vec x: = \vec x (\mathbf{U}_i) : = \mathbf{W}_1 \circ \ldots  \circ \mathbf{U}_i \circ \ldots \circ \mathbf{W}_d
\end{equation}
Therefore the original optimization problem in the large ambient space $\mathbb{R}^{n_1 \times \ldots \times n_d}$ is restricted or projected onto a relatively small subspace $\mathbb{R}^{r_{i-1} \times n_i \times r_i}$, where it can be easily solved, 
\begin{equation}
\vec W_i := \operatorname{argmin}_{\vec U_i \in \mathbb{R}^{r_{i-1} \times n_i \times r_i}} \left( \| \mathbf{W}_1 \circ \ldots  \circ \mathbf{U}_i \circ \ldots \circ \mathbf{W}_d - \mathbf{f} \|^2 \right) \ .
\end{equation}
This procedure is then continued iteratively by choosing another component $\mathbf{W}_j$ to be optimized next, resulting in a nonlinear Gau\ss{} Seidel iteration. The process of optimizing each component exactly once is often called a half-sweep.\\

Although these ideas can also be applied for the canonical format and general tensor networks as well, the tensor train and hierarchical Tucker format admit the possibility to use an orthonormal bases, which can be directly derived from the components $\vec W_i$ by corresponding (left/right) orthogonalization. With this simple post processing step, the ALS algorithm performs much better and more stable than without orthogonalization \see e.g. \cite{holtz2012alternating}. To get started, the ALS algorithm requires an initial guess, i.e\. it needs $d-1$ (left /right) orthogonal components. A usual choice is to use (Gaussian) random tensors for the $d-1$ components, possibly orthogonalized. The interesting observation is that using these random initialization one half-sweep can almost be cast to our proposed randomized TT-SVD and vice versa, in the sense that numerically exactly the same operations are performed. The only difference is that, in the picture of the randomized TT-SVD, the random tensor $\vec g$ is not chosen as a Gaussian random tensor in each step but as the first $d-i$ (contracted) random components of the ALS initialization. Note that this means that for the matrix case $d=2$ the two methods actually coincide completely. Would it be possible to extend our error bounds to the setting of using structured random tensors $\vec g$ and also to cope with the stochastic dependence implied by the fact that $\vec g$ of different steps are not sampled independently, one could for example prove stochastic error bounds for the first sweep of the ALS. Possible not only for the low rank approximation setting but also for more general objective functionals $\mathcal{J}$. Numerical results indeed do suggest that such extensions might be possible. However this is devoted to forthcoming research.

\section{Numerical Experiments}
\label{sec:numerics}

In order to provide practical proof of the performance of the presented randomized TT-SVD we conducted several numerical experiments. All calculations were performed using the \emph{xerus} C++ toolbox \cite{xerus}, which also contains our implementation of the randomized TT-SVD. The random tensors used in the following experiments are created as follows. Standard Gaussian random tensors are created by sampling each entry independently from $\mathcal{N}(0, 1)$. Sparse random tensors are created by sampling $N$ entries independently from $\mathcal{N}(0, 1)$ and placing them at positions sampled independently and evenly distributed from $[n_1] \times [n_2] \times \ldots \times [n_d]$. The low rank tensors are created by sampling the entries of the component tensors $\vec W_1, \ldots, \vec W_d$ in representation \eqref{eq:TTRepresentation}, independently from $\mathcal{N}(0, 1)$, i.e. all components $\vec W_i$ are independent standard Gaussian random tensors. In some experiments we impose a certain decay for the singular values of the involved matrifications. To this end we create a all random components as above, then for $i$ from $1$ to $d-1$ we contract $\vec W_i \circ W_{i+1}$ and then re-separate them by calculating the SVD $\vec U \vec S \vec V^T$, but replacing the $\vec S$ with a matrix $\tilde{\vec S}$ in which the singular values decay in the desired way. For a quadratical decay that is $\tilde{\vec S} := \text{diag}(1, \frac{1}{2^2}, \frac{1}{3^2}, \ldots, \frac{1}{250^2}, 0, \ldots, 0)$, where $250$ is a cut-off used in all experiments below. Then set $\vec W_i = \hat M^{-1} (\vec U)$ and $\vec W_{i+1} = \hat M^{-1}(\vec S \vec V^T)$. Note that since the later steps change the singular values of the earlier matrification, the singular values of the resulting tensor do not obey the desired decay exactly. However empirically we observed that this method yields a sufficiently well approximation for most applications, even after a single sweep.\\

A general problem is that, as described in section \ref{sec:ttFormat} the calculation of the actual best rank $\vec r$ approximation of higher order tensors is NP-hard in general. Moreover to the authors knowledge there are no non-trivial higher order tensors for which this best approximation is known in advance. Therefore a direct check of our stochastic error bound \eqref{eq:stochErrorBoundRndTTSvD} using the actual best approximation error is unfeasible. Instead most of the numerical experiments use the error of the deterministic TT-SVD introduced in section \ref{sec:ttFormat} for comparison, which gives a quasi best approximation. The factor of $\sqrt{d-1}$ is present in both error bounds, but the remaining error dependence given by \eqref{eq:stochErrorBoundRndTTSvD} is verifiable in this way.\\

If not stated otherwise we use the same values for all dimensions $n_i = n$, target ranks $r_i = r$ and (approximate) ranks of the solution $r_i^* = r^*$. In derogation from this rule the ranks are always chosen within the limits of the dimensions of the corresponding matrifications. For example for $d=8, n=4, r=20$ the actual TT-rank would be $\vec r = (4, 16, 20, 20, 20, 16, 4)$. $256$ samples are calculated for each point, unless specified otherwise. The results for the randomized TT-SVD are obtained by calculating an rank $r+p$ approximation as described in section \ref{sec:randTensorDecomp} and then use the deterministic TT-SDV to truncate this to rank $r$. This is done so that in all experiments the randomized and the deterministic approximation have identical final ranks $\vec r$.\\

\subsection{Approximation Quality for Nearly Low Rank Tensors}
\label{sec:expNearlyLowRank}

\begin{figure}
    \input{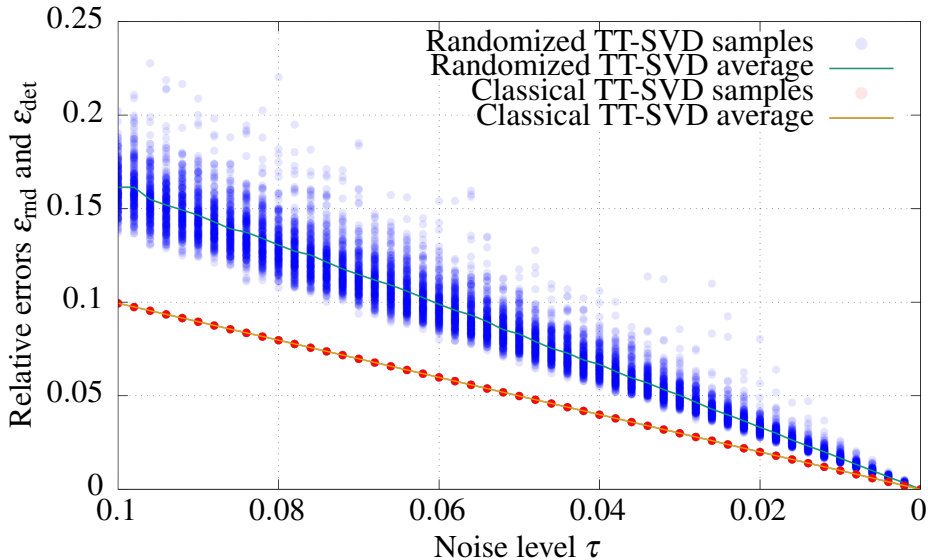}
    \caption{Approximation error of deterministic and randomized TT-SVD in dependency to the noise level. Parameters are  $d=10,\ n=4,\ r^*=10,\ r=10,\ p=5$.}
	\label{fig:noise}
\end{figure}

In this experiment we examine the approximation quality of the randomized TT-SVD for almost exact low rank tensors. I.e.\ we create random TT-rank $\vec r^*$ tensors $\vec x_{\text{exact}} \in \mathbb{R}^{n \times \ldots \times n}$ and standard Gaussian random tensors $\vec n \in \mathbb{R}^{n_1 \times \ldots \times n_d}$. The target tensor is then created as
\begin{equation}
\vec x := \frac{\vec x_{\text{exact}}}{\| \vec x_{\text{exact}} \|} + \tau \frac{\vec n}{\| \vec n \|} \ ,
\end{equation}
for some noise level $\tau$. Subsequently rank $\vec r$ approximations $\vec y_{\text{det}}$ and $\vec y_{\text{rnd}}$ of $\vec x$ are calculate using the randomized and deterministic TT-SVD. Finally we examine the relative errors
\begin{align}
\epsilon_{\text{det}} &:= \frac{\| \vec x - \vec y_{\text{det}} \|}{\| \vec x \|} & \epsilon_{\text{rnd}} &:= \frac{\| \vec x - \vec y_{\text{rnd}} \|}{\| \vec x \|} \ .
\end{align}
Figure \ref{fig:noise} show these errors for different noise levels $\tau$. The parameters are chosen as $d=10, n=4, r^*=10, r=10, p=5$, with $256$ samples calculated for each method and noise level.\\ 

As expected the error of the classical TT-SVD almost equals the noise $\tau$ for all samples, with nearly no variance. Independent of the noise level, the error $\epsilon_{\text{rnd}}$ of the randomized TT-SVD is larger by a factor of approximately $1.6$. The only exception is in the case $\tau = 0$ where both methods are exact up to numerical precision. In contrast to the classical TT-SVD, there is some variance in the error $\epsilon_{\text{rnd}}$. Notably this variance continuously decreases to zero with the noise level $\tau$. These observations are in agreement with the theoretical expectations, as theorem \ref{theo:RandTTSVDError} states that the approximation error of the randomized TT-SVD is with high probability smaller than a factor times the error of the best approximation.

\subsection{Approximation Quality with respect to oversampling}
\label{sec:expOversampling}

\begin{figure}
    \input{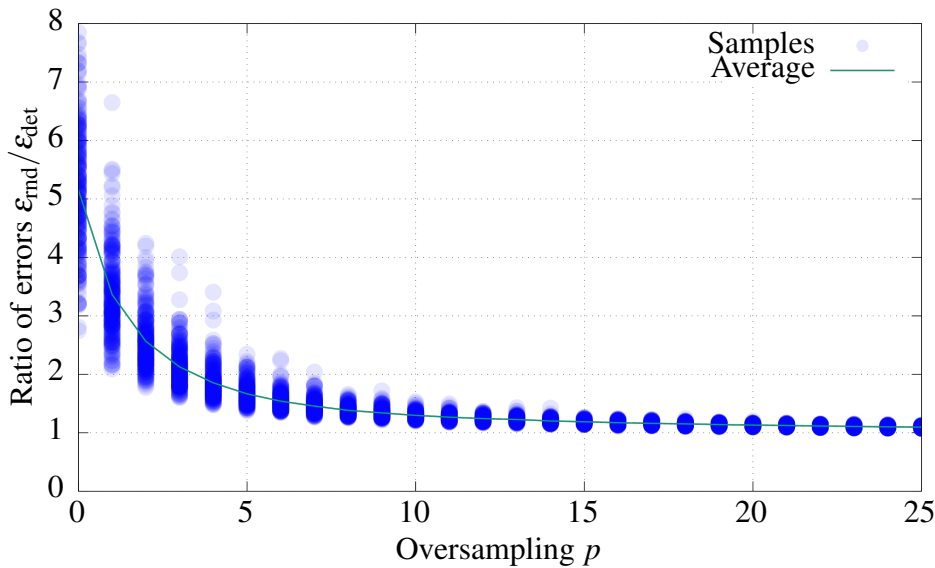}
    \caption{Approximation error of the randomized TT-SVD in terms of the deterministic error in dependency of the oversampling for nearly low rank tensors. Parameters are  $d=10,\ n=4,\ r^*=10,\ r=10,\ \tau=0.05$.}
	\label{fig:oversampling}
\end{figure}

\begin{figure}
    \input{figures/oversampling2.tex}
    \caption{Approximation error of the randomized TT-SVD in terms of the deterministic error in dependency of the oversampling for tensors with quadratically decaying singular values. Parameters are  $d=10,\ n=4,\ r=10$.}
	\label{fig:oversampling2}
\end{figure}

In this experiment we examine the influence of the oversampling parameter $p$ on the approximation quality. The first setting is similar to the one in experiment \ref{sec:expNearlyLowRank}, i.e.\ nearly low rank tensors with some noise. In contrast to experiment \ref{sec:expNearlyLowRank} the noise level $\tau = 0.05$ is fixed and the oversampling parameter $p$ is varied. For each sample we measure the error of the approximation obtained by the randomized TT-SVD $\epsilon_{\text{rnd}}$ in terms of the error of the deterministic TT-SVD $\epsilon_{\text{det}}$. The results for the parameters $d=10, n=4, r^*=10, r=10, \tau=0.05$ are shown in figure \ref{fig:oversampling}. For small $p$ a steep decent of the error factor is observed which slowly saturates towards a factor of approximately one for larger $p$. The variance decreases at a similar pace.\\

In the second setting tensors with quadratically decaying singular values are used, see the general remarks at the beginning of the section for the details of the creation. The behavior of the error factor is visualized in figure \ref{fig:degree} for the same parameters $d=10, n=4, r=10$. There are several differences compared to the first setting. Most obvious for all $p$ the factor is much smaller, i.e.\ the second setting is more favorable to the randomized TT-SVD. The same is also true for the variance. A more subtle difference, at least in the measured range of $p$, is that there are many samples for which the error factor is smaller than one, i.e.\ the randomized approximation is actually better than the deterministic one.\\

Very loosely speaking theorem \ref{theo:RandTTSVDError} predicts a $1+\frac{c}{\sqrt{p}}$ dependency of the error factor with respect to $p$, which is also roughly what is observed in both experiments.

\subsection{Approximation Quality with respect to the Order}

\begin{figure}
    \input{figures/degree.tex}
    \caption{Approximation error of the randomized TT-SVD in terms of the deterministic error in dependency of the order for nearly low rank tensors. Parameters are  $d=10,\ n=4,\ r^*=10,\ r=10,\ \tau=0.05$.}
	\label{fig:degree}
\end{figure}

\begin{figure}
    \input{figures/degree2.tex}
    \caption{Approximation error of the randomized TT-SVD in terms of the deterministic error in dependency of the order for tensors with quadratically decaying singular values. Parameters are  $d=10,\ n=4,\ r=10$.}
	\label{fig:degree2}
\end{figure}

In this third experiment the impact of the order on the approximation quality is investigated. Again the first setting uses nearly low rank tensors with some noise. The parameters are chosen as $d=10, n=4, r^*=10, r=10, \tau=0.05$. The result is shown in figure \ref{fig:degree}. As in experiment \ref{sec:expOversampling} in the second setting target tensors with quadratically decaying singular values are used. The results for the parameters $d=10, n=4, r=10$ are shown in figure \ref{fig:degree2}. For both settings the factor slightly increases from $d=4$ to $d=7$ but then stabilizes to a constant values of approximately $1.65$ and $1.95$ respectively. The same qualitative behavior is observed for the variance. This is somewhat better than expected from the theoretical results. The factor $\sqrt{d-1}$ in the error term of theorem \ref{theo:RandTTSVDError} is not visible as it is also present in the error bound of the deterministic TT-SVD. However the order also appears as an exponent in the probability, which should be observed in these results. The fact that it is not suggest that a refinement of theorem \ref{theo:RandTTSVDError} is possible in which this exponent does not appear.

\subsection{Computation time}
\begin{figure}
    \input{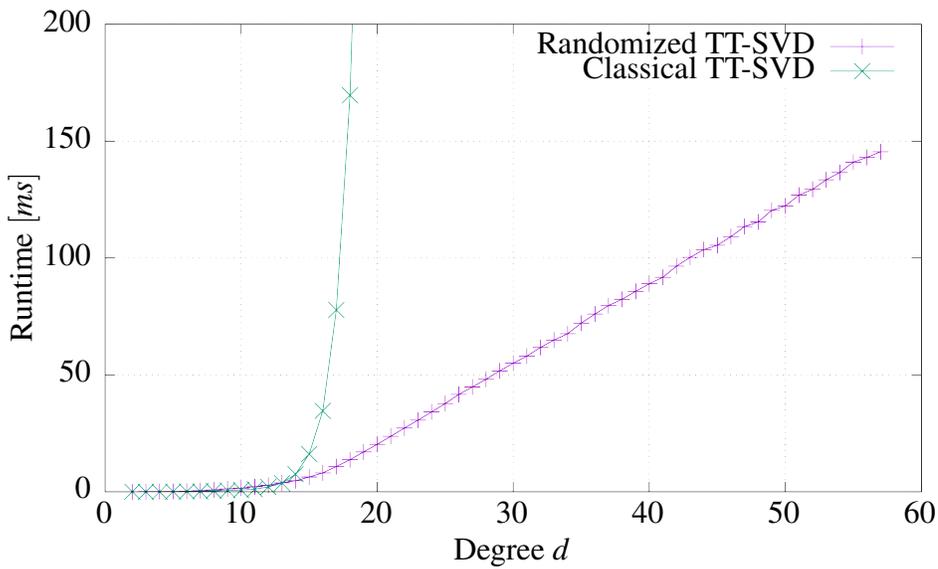}
    \caption{Run-time of the deterministic and randomized TT-SVD algorithms for different orders. Parameters are  $n=2,\ r=10,\ p=10$.}
	\label{fig:time}
\end{figure}

In this experiment we verify the computational complexity of the randomized TT-SVD algorithm, in particular the linear scaling with respect to the order for sparse tensors. To this end we create random sparse tensors with varying order and and a fixed number $N=500$ entries and measure the computation time of the TT-SVD. The other parameters are chosen as $n=2, r=10, p=10$. Figure \ref{fig:time} show the results which clearly confirm the linear scaling of the randomized TT-SVD. As a comparison also the runtime of the classical TT-SVD is given for the smaller orders. While the absolute numbers are of course hardware and implementation depended, the dramatic edge of the randomized approach is obvious.

\subsection{Approximation Quality using low rank random tensors}

\begin{figure}
    \input{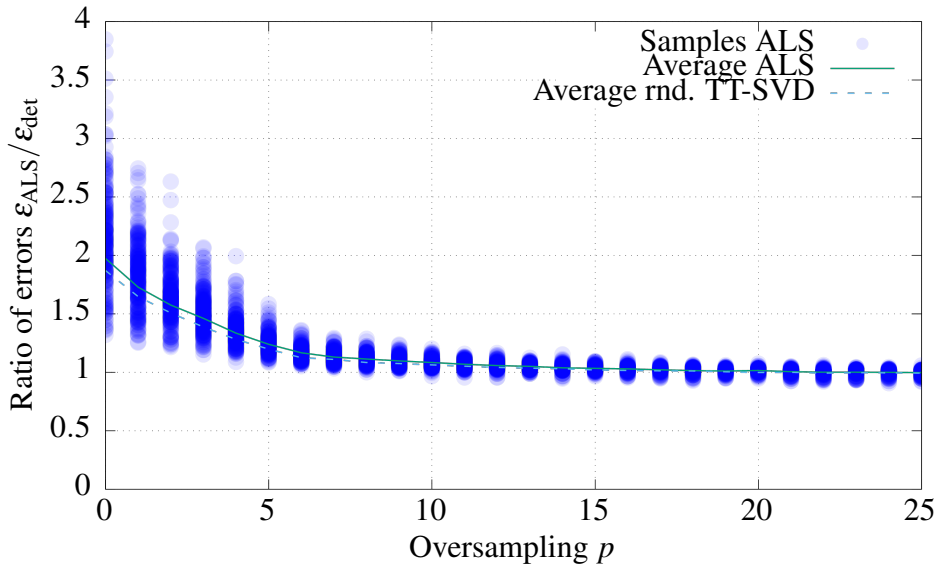}
    \caption{Approximation error of one half-sweep of the ALS with random initialization in terms of the error of the deterministic TT-SVD in dependency of the oversampling. Tensors with quadratically decaying singular values are used along the with the parameters $d=10,\ n=4,\ r=10$.}
	\label{fig:lowrank}
\end{figure}

This final experiment uses the low rank random tensor approach to the TT-SVD discussed in section \ref{sec:als}, i.e.\ instead of the proposed randomized TT-SVD one half-sweep of the ALS algorithm with random initialization is performed. The remainder of the setting is the same as the second one of experiment \ref{sec:expOversampling}, i.e.\ tensors with quadratically decaying singular values and parameters $d=10, n=4, r=10$. Figure \ref{fig:lowrank} shows the results and also as a comparison the average errors from experiment \ref{sec:expOversampling}. Apparently the error factor using the ALS half-sweep is somewhat larger than the one of the randomized TT-SVD, but otherwise exhibits the same behavior with respect to the oversampling. While there are no theoretical results on this method yet, this result is encouraging as it suggest that error bounds similar to theorem \ref{theo:RandTTSVDError} are possible for the ALS with random initialization.\\

\section{Conclusions and Outlook}
\label{sec:conclusions}

We have shown theoretically and practically that the randomized TT-SVD algorithm introduced in this work, provides a robust alternative to the classical deterministic TT-SVD algorithm at low computational expenses. In particular the randomized TT-SVD is provably exact if applied to a tensor with TT-rank smaller or equal to the target rank. For the case of actual low rank approximations stochastic error bounds hold. The numerical experiments suggest that these proven bounds are somewhat pessimistic as the observed error is mostly significantly smaller than expected. Especially we do not observe a significant deterioration of the error bound with increased order, as suggested by the current theoretical results. This leaves room for improvements and we believe that enhanced versions of our theorem are possible. On the computational side we have provided efficient implementations of the proposed algorithm, available in the \emph{xerus} C++ toolbox. For sparse tensors, we have shown that the randomized TT-SVD algorithm dramatically outperforms the deterministic algorithm, scaling only linearly instead of exponentially in the order, which was verified by measurement of the actual implementation. We believe that these results show that the randomized TT-SVD algorithm is a useful tool for low rank approximations of higher order tensors.\\

In order to avoid repetition we presented our randomized TT-SVD algorithm only for the popular tensor train format. Let us note however that the very same ideas can straight forwardly be applied to obtain an algorithm for a randomized HOSVD for the Tucker format. We expect that they can also be extended to obtain a randomized HSVD for the more general hierarchical Tucker format, but this is still work in progress. While an extension to the canonical polyadic format would certainly be desirable as well, we expect such an extension to be much more evolved, if possible at all.\\

A topic of further investigations is the use of structured random tensors in the randomized TT-SVD. For the matrix case several choices of structured random matrices are already discussed in the work of \citet{halko2011finding}. Transferring their results to the high dimensional case could allow choices of randomness which lead to reduced computational cost if the given tensor is dense, as it is the case for matrices. The even more interesting choice however is to use random low rank tensors, as already discussed in section \ref{sec:als}. On the one hand an analysis of this setting directly benefits the alternating least squares algorithm, as it would result in error bound for the first half-sweep for a random initial guess. This can be of major importance as there are mainly local convergence theories for the ALS, which is why the starting point matters a lot. On the other hand having error bounds also for this setting allow computationally fast application of the randomized TT-SVD to tensors given in various data-sparse formats, e.g.\ in the canonical, the TT or HT format and also combination of those. This is for example important for the iterative hard thresholding algorithm for tensor completion, discussed in the introduction. Here in each iteration an SVD of a low rank plus a sparse tensor has to be calculated.

\printbibliography 
\end{document}